\documentclass[12pt,a4paper]{amsart}
\usepackage[T1]{fontenc}
\usepackage[utf8]{inputenc}
\usepackage{lmodern}
\usepackage[a4paper,textwidth=148mm,hcentering,top=28mm,bottom=27mm]{geometry}
\usepackage{amsmath,amssymb,mathtools,mathrsfs}
\usepackage{enumitem}
\usepackage{tikz-cd}
\usepackage{microtype}
\usepackage{xcolor}
\usepackage[colorlinks=true,linkcolor=blue,citecolor=blue,urlcolor=blue,bookmarksnumbered=true]{hyperref}
\hypersetup{pdftitle={Algebraic integrability of p-closed corank one foliations},
 pdfauthor={Yuting Liu},
 pdfsubject={Closed rational forms, logarithmic conormal sections, transversely affine foliations, and nef anti-canonical classes}}
\setlist[enumerate,1]{label=\textup{(\roman*)},leftmargin=2.1em,itemsep=2pt,topsep=4pt}
\newtheorem{theorem}{Theorem}[section]
\newtheorem{proposition}[theorem]{Proposition}
\newtheorem{lemma}[theorem]{Lemma}
\newtheorem{corollary}[theorem]{Corollary}
\newtheorem{conjecture}[theorem]{Conjecture}

\theoremstyle{definition}
\newtheorem{definition}[theorem]{Definition}
\newtheorem{remark}[theorem]{Remark}

\theoremstyle{plain}
\newtheorem{theorema}{Theorem}

\newtheorem{theoremb}{Theorem}

\numberwithin{equation}{section}
\newcommand{\F}{\mathcal F}
\newcommand{\G}{\mathcal G}
\newcommand{\E}{\mathcal E}
\newcommand{\HH}{\mathcal H}
\newcommand{\OO}{\mathcal O}
\newcommand{\CC}{\mathbb C}
\newcommand{\QQ}{\mathbb Q}
\newcommand{\ZZ}{\mathbb Z}
\newtheorem{theoremc}{Theorem}

\newcommand{\PP}{\mathbb P}
\renewcommand{\AA}{\mathbb A}

\newcommand{\RR}{\mathbb R}

\newcommand{\Cart}{\mathcal C}
\newcommand{\sat}{\mathrm{sat}}
\newcommand{\smq}{\mathbin{\sim_{\QQ}}}
\DeclareMathOperator{\Spec}{Spec}
\DeclareMathOperator{\rank}{rank}
\DeclareMathOperator{\codim}{codim}
\DeclareMathOperator{\Sing}{Sing}
\DeclareMathOperator{\Exc}{Exc}
\DeclareMathOperator{\Supp}{Supp}
\DeclareMathOperator{\Alb}{Alb}
\DeclareMathOperator{\Pic}{Pic}

\DeclareMathOperator{\Tr}{Tr}

\DeclareMathOperator{\Frob}{Frob}

\newcommand{\NE}{\overline{\mathrm{NE}}}
\title{Algebraic integrability of $p$-closed corank one foliations}
\author{Yuting Liu}
\date{}

\begin{document}
\begin{abstract}
We verify several cases of the Ekedahl--Shepherd-Barron--Taylor
conjecture, which predicts algebraic integrability for foliations
that are $p$-closed for almost all primes.
For corank one foliations, we prove the conjecture in the following
settings: transversely affine foliations on normal complex varieties;
foliations on smooth projective varieties with nef anti-canonical
class and either canonical singularities or local holomorphic first
integrals at every point; and log canonical foliated pairs
$(\F,\Delta)$ on normal projective threefolds with
$-(K_{\F}+\Delta)$ nef.
As applications, we verify the conjecture for log canonical
degree-two foliations on $\PP^n$, $n\geq3$, log canonical corank one
Mukai foliations on smooth projective varieties of dimension at least
four and Picard number one, and Fano foliations of large index.
We also obtain an explicit criterion in the adjoint Fano setting.
\end{abstract}
\maketitle
\tableofcontents

\section{Introduction}
A basic question in the theory of foliations is whether their leaves are algebraic. A
foliation is algebraically integrable if its general leaf is open in its Zariski closure, or
equivalently if it is induced by a dominant rational map. A conjecture proposed in
\cite{ESBT99} characterizes this property in terms of reduction modulo
primes. In characteristic $p$, the $p$-th iterate of a derivation is again a derivation. A
foliation is called $p$-closed if its tangent sheaf is preserved by this operation. For a
complex foliation, the corresponding condition is imposed on the closed fibres of a
model over a finitely generated $\ZZ$-algebra; see Definition~\ref{pclosedalgdef}.
\begin{conjecture}[\cite{ESBT99}]\label{conj:esbt}
A foliation on a normal complex variety is algebraically integrable if and only if it is
$p$-closed for almost all primes.
\end{conjecture}
Conjecture~\ref{conj:esbt} remains open even for surfaces. \cite{Bos01} proved the conjecture for regular foliations defined
over number fields under the additional assumption that the general
leaves satisfy the Liouville property; more recently, \cite{Dru25}
obtained algebraicity results for separatrices of reduced surface
foliations. Conjecture~\ref{conj:esbt} was also proved for degree-one foliations on projective spaces in
\cite[Lemma~11.5]{MP24}. For corank one foliations with
numerically trivial canonical class and canonical singularities on projective klt varieties, the conjecture
follows from \cite{Dru21} and \cite[Section~8]{DO22}.
For other related
results, see \cite{And04,BCL09,LL25,LL26,Liu26}.

In this paper, we prove algebraic integrability for transversely affine
corank one foliations on normal complex varieties which are $p$-closed
for almost all primes. We also prove the conjecture on smooth
projective varieties when $-K_{\F}$ is nef and $\F$ either has canonical singularities
or admits local holomorphic first integrals at every point. For log canonical
foliated pairs $(\F,\Delta)$ on normal projective threefolds, we prove the
conjecture when $-(K_{\F}+\Delta)$ is nef. Finally, we consider Fano foliations,
degree-two foliations on projective space, and adjoint Fano foliated
structures.

\begin{theorema}\label{main:forms}
Let $X$ be a normal complex variety, and let $\F$ be a transversely
affine corank one foliation on $X$. If $\F$ is $p$-closed for almost
all primes, then $\F$ is algebraically integrable.
\end{theorema}

In particular, the theorem applies to foliations defined by closed
rational one-forms.
Theorem~1.5 of \cite{Gra26} considers vector fields
$f(x)\partial_x+g(y)\partial_y$ with $f,g\in\QQ(t)$.
In dimension two these fields are annihilated by the closed rational
one-form $dx/f(x)-dy/g(y)$, so the corresponding foliations are
covered by Theorem~\ref{main:forms}.
The higher-dimensional results of \cite{Gra26} consider rank-one
rather than corank one foliations.

A corank one foliation $\F$ is transversely affine if there exist
rational one-forms $\omega\neq0$ and $\eta$ such that
\[
\F=\ker\omega,\qquad
d\omega=\omega\wedge\eta,\qquad d\eta=0;
\]
see Definition~\ref{def:affine}.
On a smooth projective variety, this is equivalent to the existence
of a meromorphic flat connection on $N_{\F}$ compatible with the
defining twisted one-form.
Away from the singular and polar loci, local first integrals can
be chosen with affine transition maps $z\mapsto az+b$, where
$a\in\CC^*$ and $b\in\CC$ are constant.

Our second result is to prove Conjecture \ref{conj:esbt} for foliations with nef anti-canonical divisor under mild singularity conditions.
\begin{theoremb}\label{main:nef}
Let $X$ be a smooth complex projective variety of dimension at least two, and let
$\F\subset T_X$ be a corank one foliation. Assume that $-K_{\F}$ is nef, that
$\F$ is $p$-closed for almost all primes, and that at least one of the following holds:
\begin{enumerate}
\item $\F$ has canonical singularities;
\item $\F$ admits a local holomorphic first integral at every point.
\end{enumerate}
Then $\F$ is algebraically integrable and regular. Moreover, there is a regular
rank-one foliation $\E$ such that
\[
T_X=\F\oplus\E.
\]
\end{theoremb}
In particular, by \cite[Corollary~1.7]{Liu26}, we obtain algebraic integrability
for corank one foliations with nef anti-canonical class and only simple singularities,
provided that they are $p$-closed for almost all primes. The same conclusion holds
when the singular locus has codimension at least three. In both cases the foliation
is in fact regular; see Corollary~\ref{cor:small-sing}.

In dimension three, we prove the following result for log canonical
foliated pairs.

\begin{theoremc}\label{main:log-nef}
Let $X$ be a normal complex projective threefold, and let $\F$ be
a corank one foliation on $X$ which is $p$-closed for almost all primes.
Let $\Delta$ be an effective $\QQ$-divisor such that $(\F,\Delta)$
is log canonical. If $-(K_{\F}+\Delta)$ is nef, then $\F$ is
algebraically integrable.
\end{theoremc}

\par\medskip
\noindent\textbf{Applications.}
We prove Conjecture~\ref{conj:esbt} for several classes of Fano
foliations.
Recall that a foliation $\F$ is $\QQ$-Fano if $K_{\F}$ is
$\QQ$-Cartier and $-K_{\F}$ is ample.
Its index $\iota(\F)$ is the largest positive rational number $a$
such that $-K_{\F}\smq aH$ for an ample Cartier divisor $H$.

\begin{corollary}\label{cor:intro-fano}
Let $X$ be a $\QQ$-factorial klt complex projective variety of
dimension $n\geq3$, and let $\F\subset T_X$ be a corank one
$\QQ$-Fano foliation which is $p$-closed for almost all primes.
Then $\F$ is algebraically integrable if either of the following holds:
\begin{enumerate}
\item $\iota(\F)\geq n-2$;
\item $\F$ has log canonical singularities, and either $n=3$ or $X$ is
a smooth fourfold with $\rho(X)=1$.
\end{enumerate}
\end{corollary}

In dimension three, the index assumption in (i) is automatic
when $K_{\F}$ is Cartier.

We also prove the conjecture for log canonical degree-two
foliations and log canonical Mukai foliations.
Recall that a Fano foliation of rank $r\geq3$ is called a
Mukai foliation if its index is $r-2$.

\begin{corollary}\label{cor:intro-degree-mukai}
Let $X$ be a smooth complex projective variety of dimension
$n\geq3$, and let $\F\subset T_X$ be a corank one foliation
with log canonical singularities.
Assume that $\F$ is $p$-closed for almost all primes.
Then $\F$ is algebraically integrable if either of the following holds:
\begin{enumerate}
\item $X\simeq\PP^n$ and $\F$ has degree two;
\item $n\geq4$, $\rho(X)=1$, and $\F$ is a Mukai foliation.
\end{enumerate}
\end{corollary}

Theorem~\ref{main:nef} also gives an algebraic integrability
criterion in the adjoint Fano setting.

\begin{corollary}\label{cor:intro-adjoint}
Let $X$ be a smooth complex projective variety of dimension
$n\geq2$, and let $\F\subset T_X$ be a corank one foliation
which is $p$-closed for almost all primes.
Assume that $\F$ either has canonical singularities or admits
a local holomorphic first integral at every point.
If there exists a rational number $t$ such that
\[
\frac{n+1}{n+2}\leq t<1
\qquad\text{and}\qquad
-\bigl(tK_{\F}+(1-t)K_X\bigr)\ \text{is ample},
\]
then $\F$ is algebraically integrable and regular.
\end{corollary}

\par\medskip \noindent \textbf{Outline of proofs.}
To prove Theorem~\ref{main:forms}, we first consider foliations defined
by closed rational one-forms in Theorem~\ref{thm:closed}.
We pass to a smooth projective model on which the polar divisor
of the defining closed rational one-form $\omega$ has simple
normal crossing support.
We consider the vector space spanned by $\omega$ and the global
logarithmic one-forms with poles along this support.
If the foliation has no nonconstant rational first integral,
the kernel of exterior product with $\omega$ on this space
is one-dimensional.

After reduction modulo almost all primes, this kernel remains
one-dimensional and the Cartier operator sends the reduced form to a logarithmic
form. The $p$-closedness assumption therefore forces the image of the Cartier operator
to be a constant multiple of the reduced form.
If $\omega$ has a pole of order greater than one, this constant
vanishes. On the complement of the polar divisor, we use $\omega$
to equip the trivial rank-two bundle with a flat connection that is
an extension of the trivial connection by itself.
The Cartier identity \cite[Lemma~4.16]{MP24} shows that its reductions
have zero $p$-curvature. Andr\'e's extension theorem
\cite[Corollaries~4.3.4 and~7.1.3]{And04} then yields the existence of a rational
first integral.

If $\omega$ is logarithmic with a nonzero residue, we divide $\omega$
by that residue. The resulting form has reductions fixed by the Cartier operator
and defines a flat connection on the trivial line bundle over the
complement of the polar divisor.
The same Cartier identity shows $p$-curvature vanishes, so the
rank-one Grothendieck--Katz theorem
\cite[Corollaries~4.3.6 and~7.1.3]{And04} yields the existence of a rational first integral.

Finally, if $\omega$ is holomorphic, the Albanese morphism reduces
the problem to a linear foliation on an abelian variety.
Its algebraic integrability follows from
\cite[Proposition~10.7]{Dru21}.

For the general case of Theorem~\ref{main:forms}, the rational forms
defining the transverse affine structure allow us to pass to a
smooth projective model.
The structure theorem \cite[Theorem~2.2]{CLNP22} then reduces the
problem to a foliation defined by a closed rational one-form
after a generically finite cover or to a transversely affine
Riccati foliation on a ruled surface. In the Riccati case, we first show that the polar divisor of the transverse affine
connection has a component dominating the base curve which is an invariant multisection.
After a finite base change, it gives an invariant
section $\Sigma$. Over a suitable open
subset of the base curve, pushing forward $\OO(\Sigma)$ along the ruling defines
a rank-two vector bundle. Horizontal differentiation defines a flat connection
on this bundle whose projectivization recovers the Riccati foliation.
The constant functions define a trivial rank-one subconnection.
The $p$-closedness assumption implies that the rank-two connection, and hence
its rank-one quotient, has zero $p$-curvature after reduction modulo almost
all primes. The rank-one Grothendieck--Katz theorem
\cite[Corollaries~4.3.6 and~7.1.3]{And04}
shows the quotient connection is isotrivial. Andr\'e's extension theorem
\cite[Corollaries~4.3.4 and~7.1.3]{And04}
then implies that the rank-two connection is isotrivial.
After a finite \'etale base change, its projectivization is the trivial
$\PP^1$-bundle with the product connection. Its horizontal leaves are
therefore algebraic, and algebraic integrability descends to the original
foliation.

For Theorem~\ref{main:nef}, let $\HH=\F^a$ be the algebraic part.
If $\HH=0$, then \cite[Corollary~4.13]{CP19} and nefness of
$-K_{\F}$ imply $K_{\F}\equiv0$.
Since $\F$ is not uniruled, \cite[Corollary~3.8]{LPT18} implies
that it is canonical. Theorem~\ref{thm:numtriv} then implies
algebraic integrability, contradicting $\HH=0$.

Thus $\HH$ has positive rank. It has a compact leaf by
Lemma~\ref{lem:canonical-compact} under hypothesis (i), and by
Lemma~\ref{lem:local-compact} under hypothesis (ii).
Proposition~\ref{prop:comparison} shows $K_{\HH}\equiv K_{\F}$.
Applying Theorem~\ref{thm:ou} to $\HH$ and descending $\F$ along
the resulting smooth fibration $q:X\to Y$, we obtain
$\F=q^{-1}\G$ with $K_{\G}\equiv0$.
Moreover, the algebraic part of $\G$ is regular and has
numerically trivial canonical class.

By \cite[Corollary~3.8]{LPT18}, this algebraic part is not uniruled,
so neither is $\G$.
If $\G\neq0$, the same result implies that $\G$ is canonical.
$p$-closedness descends to $\G$, and
Theorem~\ref{thm:numtriv} proves its algebraic integrability.
The case $\G=0$ is immediate.
Thus $\F$ is algebraically integrable and equals $\HH$.
In particular, $\F$ has a compact leaf, so
\cite[Corollary~1.2]{Ou25} implies regularity and the involutive
splitting of $T_X$.

For Theorem~\ref{main:log-nef}, we first take an F-terminalisation
\cite[Theorem~11.1]{CS21} which is a $\QQ$-factorial klt model $Y$ on
which $\F_Y$ is canonical and $-(K_{\F_Y}+B)$ is nef for an effective
boundary $B$. If the algebraic part $\HH=\F_Y^a$ has rank zero,
\cite[Corollary~4.13]{CP19} reduces the problem to
Theorem~\ref{thm:numtriv}. Otherwise, if $\F_Y$ is not algebraically
integrable, $\HH$ has rank one and is induced by an almost holomorphic
map to a surface by \cite[Theorem~5.4]{CS25}. 
Using \cite[Proposition~4.1]{Dru17}, we prove in
Lemma~\ref{lem:pair-algebraic-comparison} that
$K_{\HH}-K_{\F_Y}$ is pseudo-effective. We take a minimal model $(S,\G)$ of the quotient foliation,
where $S$ is klt, $\G$ has canonical singularities, and
$K_{\G}$ is nef. Then the negativity
lemma together with pseudo-effectivity of $K_{\HH}-K_{\F_Y}$ imply $K_{\G}\equiv0$; see
Lemma~\ref{lem:pair-nef-surface-quotient}. Theorem~\ref{thm:numtriv}
shows $\G$ algebraically integrable, contradicting pure transcendence.

For Corollary~\ref{cor:intro-degree-mukai}(i), the classification in
\cite[Theorem~A]{CLN96} shows that either the foliation is defined
by a closed rational one-form, or it is the linear pull-back of
a degree-two foliation on $\PP^2$. We exclude the second possibility
using log canonicity: blowing up the centre of the projection produces
a non-invariant exceptional divisor with discrepancy $-2$; see
Lemma~\ref{lem:projection}. Theorem~\ref{main:forms} then applies to the remaining case.

For Corollary~\ref{cor:intro-fano}(i), suppose that $\F$ is not
algebraically integrable. Its algebraic rank is then at most $n-2$,
whereas the index bound \cite[Theorem~1.3]{Liu23} implies
$r^a(\F)\geq\iota(\F)\geq n-2$.
Thus equality holds. The equality case of the same
theorem shows that $X$ is a normal generalized cone over a surface
$T$ and that $\F$ is the pull-back of a foliation $\G$ on $T$
with $K_{\G}\smq0$.
We show that $T$ is klt and that $\G$ is purely transcendental and
$p$-closed for almost all primes. Since $\G$ is purely transcendental
and $K_{\G}\equiv0$, it has canonical singularities.
The algebraic integrability results in \cite{Dru21} and
\cite[Section~8]{DO22} then imply that $\G$ is algebraically
integrable, contradicting pure transcendence.

For Corollary~\ref{cor:intro-adjoint}, the main step is to prove
nefness of $-K_{\F}$ from ampleness of $-(tK_{\F}+(1-t)K_X)$.
If $-K_{\F}$ were not nef, we could find a $K_X$-negative extremal ray
on which $K_{\F}$ is positive. The ordinary cone theorem provides a rational curve $C$ spanning this ray
with $0<-K_X\cdot C\leq n+1$. Since $K_{\F}$ is Cartier,
$K_{\F}\cdot C\geq1$. For $t\geq(n+1)/(n+2)$, these inequalities contradict
the assumed ampleness. When the adjoint divisor is anti-ample for all
$t$ sufficiently close to one, nefness follows directly by taking a limit.
In either case, Theorem~\ref{main:nef} applies.

\noindent \textbf{Acknowledgments.}
The author thanks his PhD supervisor, Paolo Cascini, for
his constant support and many valuable suggestions, and Calum Spicer for many
fruitful discussions. He is also grateful to Federico Bongiorno, St\'ephane Druel, Jie
Liu, Roktim Mascharak, and Jorge Vit\'orio Pereira for helpful discussions. The
author is supported by Imperial College London and the Heilbronn Institute for
Mathematical Research.

\noindent \textbf{AI disclosure}. All mathematical ideas in this paper are due to the author. The author used ChatGPT 6 Astra to fix a gap in the original proof of Theorem \ref{main:forms}. The author also used ChatGPT to check all the arguments and improve the exposition of this paper. The author takes full responsibility for the contents of this paper.
\section{Preliminaries}
Unless otherwise stated, varieties are integral complex algebraic varieties. All
varieties to which a global positivity assumption is applied are projective. We
distinguish linear equivalence, $\QQ$-linear equivalence, and numerical equivalence by
$\sim$, $\smq$, and $\equiv$.

\subsection{Foliations and singularities}
\begin{definition}
Let $X$ be a normal variety. A foliation $\F$ of rank $r$ on $X$ is a coherent
subsheaf $\F\subset T_X$ of rank $r$ satisfying:
\begin{enumerate}
\item $\F$ is closed under the Lie bracket,
\item $\F$ is saturated in $T_X$, that is, $T_X/\F$ is torsion-free.
\end{enumerate}
The corank of $\F$ is $\dim X-r$.
\end{definition}
\begin{remark}
A foliation on a normal variety is uniquely determined by its restriction to any
nonempty Zariski open subset.
\end{remark}
\begin{definition}
Let $\F$ be a foliation of rank $r$ on a normal variety $X$. Its canonical
divisor $K_{\F}$ and normal sheaf $N_{\F}$ are defined by
\begin{equation}\label{eq:canonical}
\OO_X(-K_{\F})\simeq\det\F:=\left(\bigwedge^r\F\right)^{**},
\qquad N_{\F}:=(T_X/\F)^{**}.
\end{equation}
The dual sheaf $N_{\F}^*$ is called the conormal sheaf of $\F$.
\end{definition}
\begin{definition}
Let $X$ be smooth and let $\F$ be a foliation on $X$. The foliation is regular at
$x\in X$ if $\F$ is a subbundle of $T_X$ in a neighbourhood of $x$. Its singular
locus, denoted by $\Sing(\F)$, is the complement of its regular locus.
\end{definition}
\begin{remark}\label{omegaF}
Let $\F$ be a corank one foliation on a normal variety $X$. The natural morphism
$N_{\F}^*\to\Omega_X^{[1]}$, where $\Omega_X^{[1]}:=(\Omega_X^1)^{**}$, determines
a nonzero twisted reflexive one-form
\[
\omega_{\F}\in H^0\bigl(X,(\Omega_X^{[1]}\otimes N_{\F})^{**}\bigr)
\]
with no divisorial zeroes. On $X_{\mathrm{reg}}$, a local trivialization of $N_{\F}$
identifies $\omega_{\F}$ with a holomorphic one-form satisfying
$\omega_{\F}\wedge d\omega_{\F}=0$.

Conversely, a nonzero rational one-form $\omega$ satisfying $\omega\wedge d\omega=0$
defines a corank one foliation by taking the saturated kernel of contraction with
$\omega$; see \cite[Sections~3.1--3.2]{Dru21}.
\end{remark}
\begin{definition}
Let $\varphi:X\dashrightarrow Y$ be a dominant rational map between normal
varieties, and let $\G$ be a foliation on $Y$. On suitable dense smooth open subsets
on which $\varphi$ is a smooth morphism, the inverse image of $\G$ under
$d\varphi$ defines a foliation. Its saturated extension to $X$ is called the
pull-back foliation and is denoted by $\varphi^{-1}\G$. When $\varphi$ is birational,
it is also called the strict transform of $\G$.
\end{definition}
\begin{definition}
Let $X$ be a normal variety and let $\F$ be a foliation on $X$. A subvariety
$Z\subset X$ is $\F$-invariant if every local section $v$ of $\F$ preserves its
ideal sheaf, that is, $v(I_Z)\subset I_Z$.
\end{definition}
\begin{definition}
Let $\F$ be a foliation on a normal variety $X$ with $K_{\F}$ $\QQ$-Cartier.
Let $\pi:X'\to X$ be a projective birational morphism from a normal variety,
and set $\F'=\pi^{-1}\F$. Choosing compatible canonical divisors, write
\begin{equation}\label{eq:discrepancies}
K_{\F'}=\pi^*K_{\F}+\sum_E a(E,\F)E,
\end{equation}
where the sum runs over the $\pi$-exceptional prime divisors. The coefficient
$a(E,\F)$ is the discrepancy of $\F$ along $E$.

For a prime divisor $E$ over $X$, put
\[
\epsilon_{\F}(E)=
\begin{cases}
0,&\text{if $E$ is invariant under the induced foliation,}\\
1,&\text{otherwise.}
\end{cases}
\]
The foliation $\F$ has canonical, respectively log canonical, singularities if
\[
a(E,\F)\geq0,\qquad\text{respectively}\qquad
a(E,\F)\geq-\epsilon_{\F}(E),
\]
for every prime divisor $E$ exceptional over $X$.

More generally, a foliated pair $(\F,\Delta)$ consists of a foliation
$\F$ on a normal variety $X$ and an effective $\QQ$-divisor $\Delta$
such that $K_{\F}+\Delta$ is $\QQ$-Cartier. Its discrepancies are
defined by
\[
K_{\F'}+\pi_*^{-1}\Delta
=\pi^*(K_{\F}+\Delta)+\sum_E a(E,\F,\Delta)E,
\]
where the sum runs over the $\pi$-exceptional prime divisors. For a
prime divisor $D$ on $X$, put $a(D,\F,\Delta)=-\operatorname{coeff}_D\Delta$.
The pair $(\F,\Delta)$ is log canonical if
\[
a(E,\F,\Delta)\geq-\epsilon_{\F}(E)
\]
for every prime divisor $E$ over $X$; see \cite[Section~2.4]{CS21}.
\end{definition}
\begin{definition}
Let $X$ be a smooth complex variety, let $\F$ be a corank one foliation on $X$,
and let $x\in X$. A local holomorphic first integral of $\F$ at $x$ is a
nonconstant holomorphic function $h:U\to\CC$ on a connected analytic neighbourhood
$U\subset X^{\mathrm{an}}$ of $x$ such that
$\F^{\mathrm{an}}|_U=(\ker dh)^{\sat}$.
\end{definition}

\subsection{Algebraic parts and compact leaves}
\begin{definition}
Let $X$ be a normal projective variety and let $\F$ be a foliation on $X$.
The foliation $\F$ is algebraically integrable if its general leaf is algebraic,
or equivalently if it is induced by a dominant rational map.

The foliation $\F$ is purely transcendental if no positive-dimensional algebraic
subvariety through a general point of $X$ is tangent to $\F$.
\end{definition}
\begin{proposition}[{\cite[Definition~2.3]{Dru17}; \cite[Definition~3.2]{Liu23}}]\label{prop:algebraic-part}
Let $X$ be a normal projective variety and let $\F$ be a foliation on $X$.
There exist a normal projective variety $Z$, a dominant rational map
$\varphi:X\dashrightarrow Z$ with connected general fibres, and a purely
transcendental foliation $\G$ on $Z$ such that
\begin{equation}\label{eq:quotient}
\F=\varphi^{-1}\G.
\end{equation}
The variety $Z$ is unique up to birational equivalence.
\end{proposition}
\begin{definition}
In the setting of Proposition~\ref{prop:algebraic-part}, the foliation induced by
$\varphi$ is called the algebraic part of $\F$ and is denoted by $\F^a$.
Its rank $r^a(\F):=\rank\F^a$ is called the algebraic rank of $\F$.
The foliation $\G$ is called the transcendental part of $\F$.
\end{definition}
\begin{proposition}[{\cite[Corollary~4.13]{CP19}}]
Let $X$ be a smooth projective variety and let $\G$ be a purely transcendental
foliation on $X$. Then $K_{\G}$ is pseudo-effective.
\end{proposition}
\begin{definition}
A dominant rational map $\varphi:X\dashrightarrow Z$ is almost holomorphic if
there exist dense Zariski open subsets $X^\circ\subset X$ and $Z^\circ\subset Z$
such that $\varphi|_{X^\circ}:X^\circ\to Z^\circ$ is a proper surjective morphism.

A compact leaf of a foliation on a smooth complex variety is a leaf contained in
the regular locus of the foliation which is compact in the analytic topology.
\end{definition}
\begin{proposition}[{\cite[Theorem~5.4]{CS25}}]\label{prop:almost-holo}
Let $X$ be a $\QQ$-factorial klt projective variety and let $\F$ be a foliation with
canonical singularities. Then the algebraic part $\F^a$ is induced by an almost
holomorphic map.
\end{proposition}
\begin{lemma}\label{lem:canonical-compact}
Let $X$ be a smooth projective variety and let $\F$ be a foliation with canonical
singularities. If $r^a(\F)>0$, then $\F^a$ has a compact leaf.
\end{lemma}
\begin{proof}
By Proposition~\ref{prop:almost-holo}, there is an almost holomorphic map
$\varphi:X\dashrightarrow Z$ inducing $\F^a$, with connected general fibres.
Choose dense open subsets $X^\circ\subset X$ and $Z^\circ\subset Z$ such that
$\varphi^\circ:X^\circ\to Z^\circ$ is a proper surjective morphism.
After shrinking $Z^\circ$ and replacing $X^\circ$ by its inverse image, generic
smoothness and properness allow us to assume that $Z^\circ$ is smooth and
$\varphi^\circ$ is smooth with connected fibres. On this open subset,
$\F^a|_{X^\circ}=T_{X^\circ/Z^\circ}$. Every fibre of $\varphi^\circ$ is therefore
a smooth connected projective leaf of $\F^a$.
\end{proof}

\subsection{\texorpdfstring{$p$}{p}-closeness for almost all primes}
We follow the terminology of \cite[Section~3.1]{Liu26}.

\begin{definition}
Let $k$ be a field of characteristic $p>0$, let $X$ be a normal variety
over $k$, and let $\F$ be a foliation on $X$. For a local vector field $v$,
write $v^{[p]}$ for its $p$-fold composition as a derivation; in
characteristic $p$, this is again a derivation.
The foliation $\F$ is \emph{$p$-closed} if
$v^{[p]}\in H^0(U,\F)$ for every open subset $U\subset X$ and every
$v\in H^0(U,\F)$.
\end{definition}

\begin{definition}
Let $X$ be a normal variety over a field of characteristic $p>0$,
and let $\F$ be a foliation on $X$. The $p$-th power operation induces
an $\OO_X$-linear morphism
\[
\psi_{\F}:\Frob_X^*\F\longrightarrow (T_X/\F)^{**},
\qquad
1\otimes v\longmapsto v^{[p]}\bmod\F,
\]
where $\Frob_X$ denotes the absolute Frobenius morphism.
This morphism is called the \emph{$p$-curvature} of $\F$.
The foliation $\F$ is $p$-closed if and only if $\psi_{\F}=0$.
\end{definition}

\begin{definition}\label{pclosedalgdef}
Let $K$ be a field finitely generated over $\QQ$, let $X$ be a normal
variety over $K$, and let $\F$ be a foliation on $X$.
We say that $\F$ is \emph{$p$-closed for almost all primes} if there exist:
\begin{enumerate}[label=\textup{(\arabic*)}]
\item a finitely generated integral $\ZZ$-subalgebra $R\subset K$
with $\operatorname{Frac}(R)=K$;
\item a flat morphism of finite type
\[
\mathscr X\longrightarrow S:=\Spec R
\]
with normal fibres, together with an identification of its generic
fibre with $X$;
\item a coherent subsheaf
\[
\mathscr F\subset T_{\mathscr X/S},
\]
flat over $S$, whose generic fibre is identified with $\F$ and whose
fibres are foliations of rank $\rank\F$;
\end{enumerate}
such that $\mathscr F_s$ is $p(s)$-closed for every closed point $s$
of some nonempty open subset $S^\circ\subset S$.
Here $k(s)$ is the residue field of $s$ and
$p(s)=\operatorname{char}k(s)$.

For a foliation $\F$ on a normal complex variety $X$, we require
the above condition for a pair $(X_K,\F_K)$ defined over some
finitely generated subfield $K\subset\CC$ whose base change to
$\CC$ is $(X,\F)$.

We call the full model
\[
(\mathscr X,\mathscr F)\longrightarrow S
\]
a \emph{spreading out} of $(X,\F)$.
\end{definition}

\begin{remark}
By enlarging the finitely generated ring of definition and shrinking
the base, we may include finitely many additional algebraic data in
a spreading out. These operations preserve $p$-closedness for almost
all primes, since vanishing of $p$-curvature is preserved by extension
of the residue field.
\end{remark}
\begin{lemma}[{\cite[Lemmas~3.3, 3.8 and~3.24]{Liu26}}]\label{lem:open}
Let $X$ be a normal variety over a field of characteristic $p>0$,
let $\F$ be a foliation on $X$, and let $U\subset X$ be a nonempty
open subset. Then $\F$ is $p$-closed if and only if $\F|_U$ is
$p$-closed.

For a foliation on a normal complex variety, the property of being
$p$-closed for almost all primes can be checked on any
nonempty Zariski open subset. In particular, it is invariant under
birational transform.
\end{lemma}

\begin{lemma}\label{lem:dominant}
Let $q:X\dashrightarrow Y$ be a dominant rational map between normal complex
varieties, let $\G$ be a foliation on $Y$, and set $\F=q^{-1}\G$.
\begin{enumerate}
\item If $\F$ is $p$-closed for almost all primes, then so is $\G$.
\item If $\G$ is algebraically integrable, then so is $\F$.
\end{enumerate}
\end{lemma}
\begin{proof}
For (i), restrict to smooth dense open subsets on which $q$ is a smooth
surjective morphism and $\G$ is a subbundle of $T_Y$. Choose a model of these
data. After shrinking its base, the same properties hold on every closed fibre.

Let $s$ be such a closed point, put $p=p(s)$, and let $v$ be a local section of
$\G_s$. Smoothness implies, locally on $X_s$, a $q_s$-related vector field
$\widetilde v$ belonging to $\F_s$. The vector field $\widetilde v^{[p]}$ is then
$q_s$-related to $v^{[p]}$. Since $\F_s$ is $p$-closed, $\widetilde v^{[p]}$
belongs to $\F_s$, so the image of $v^{[p]}$ in $T_{Y_s}/\G_s$ vanishes after
pull-back by $q_s$. Faithful flatness implies that this image is zero. Thus $\G_s$
is $p$-closed, and Lemma~\ref{lem:open} shows the assertion on $Y$.

For (ii), let $h:Y\dashrightarrow Z$ be a dominant rational map inducing $\G$.
Then $h\circ q$ induces $\F$.
\end{proof}
\begin{lemma}\label{lem:finite}
Let $\pi:Y\dashrightarrow X$ be a dominant generically finite rational map
between normal complex varieties. Let $\F$ be a foliation on $X$ and set
$\F_Y=\pi^{-1}\F$. Then:
\begin{enumerate}
\item $\F$ is $p$-closed for almost all primes if and only if $\F_Y$ is
$p$-closed for almost all primes;
\item $\F$ is algebraically integrable if and only if $\F_Y$ is algebraically integrable.
\end{enumerate}
\end{lemma}
\begin{proof}
Choose smooth dense open subsets on which $\pi$ is finite and etale. After
spreading out and shrinking the base, the induced morphisms on closed fibres are
also finite and etale. On these open subsets, the $p$-curvature morphism of the
pull-back foliation is the pull-back of the original $p$-curvature morphism.
Faithful flatness and Lemma~\ref{lem:open} prove (i).

Algebraic integrability pulls back by Lemma~\ref{lem:dominant}. Conversely,
suppose that $\F_Y$ is algebraically integrable. The closure of the image of a
general algebraic leaf of $\F_Y$ has dimension $\rank\F$ and is tangent to $\F$
at its general point. These images cover a dense open subset of $X$, so $\F$
is algebraically integrable.
\end{proof}

\subsection{The Cartier operator}
Let $k$ be a perfect field of characteristic $p>0$, and let $X$ be a smooth
$k$-variety. Denote by
\[
Z\Omega^1_{X/k}:=\ker(d:\Omega^1_{X/k}\to\Omega^2_{X/k}),
\qquad
B\Omega^1_{X/k}:=\operatorname{im}(d:\OO_X\to\Omega^1_{X/k})
\]
the sheaves of closed and locally exact one-forms, respectively. Their direct
images under the absolute Frobenius morphism $\Frob_X$ are $\OO_X$-modules,
with scalar multiplication $a\cdot\eta=a^p\eta$.
\begin{proposition}[Cartier; {\cite[Section~4.4]{MP24}}]\label{prop:cartier}
There is a unique surjective $\OO_X$-linear morphism
\[
\Cart_X:(\Frob_X)_*Z\Omega^1_{X/k}\longrightarrow\Omega^1_{X/k}
\]
such that
\[
\ker\Cart_X=(\Frob_X)_*B\Omega^1_{X/k},
\qquad \Cart_X(f^{p-1}df)=df
\]
for every local function $f$. It induces an isomorphism
\[
\frac{(\Frob_X)_*Z\Omega^1_{X/k}}{(\Frob_X)_*B\Omega^1_{X/k}}
\xrightarrow{\ \sim\ }\Omega^1_{X/k}.
\]
\end{proposition}
The morphism $\Cart_X$ is called the absolute Cartier operator. Viewed as an
additive operator on closed one-forms, it satisfies
\[
\Cart_X(a^p\eta)=a\,\Cart_X(\eta).
\]
It also acts on closed rational one-forms by restriction to the generic point.
We write $\Cart$ when the underlying variety is clear.
\begin{corollary}\label{cor:cartier-rules}
Let $h\in k(X)$, $f\in k(X)^\times$, and $\lambda\in k$. Then
\begin{equation}\label{eq:cartier-rules}
\Cart(dh)=0,\qquad
\Cart\left(\lambda\frac{df}{f}\right)=\lambda^{1/p}\frac{df}{f}.
\end{equation}
\end{corollary}
\begin{proof}
The first identity follows from the description of the kernel. For the second, write
\[
\lambda\frac{df}{f}
=\left(\frac{\lambda^{1/p}}{f}\right)^p f^{p-1}df
\]
and apply Proposition~\ref{prop:cartier} and Frobenius semilinearity.
\end{proof}
\begin{lemma}\label{lem:cartier-functor}
Let $f:Y\to X$ be a morphism between smooth $k$-varieties. For every closed
regular one-form $\eta$ on an open subset of $X$, one has
\[
\Cart_Y(f^*\eta)=f^*\Cart_X(\eta)
\]
on its inverse image. The same identity holds for rational forms whenever their
pull-backs are defined.
\end{lemma}
\begin{proof}
The inverse Cartier isomorphism sends $dg$ to the de Rham cohomology
class of $g^{p-1}dg$. Since pull-back commutes with multiplication
and the exterior differential, it commutes with the inverse Cartier
isomorphism, and hence also with the Cartier isomorphism.
\end{proof}
\begin{lemma}[{\cite[Lemma~4.16]{MP24}}]\label{lem:cartier-contraction}
Let $\eta$ be a closed rational one-form on $X$, and let $v$ be a rational vector
field. Then
\begin{equation}\label{eq:cartier-contraction}
\bigl(\iota_v\Cart_X(\eta)\bigr)^p
=\iota_{v^{[p]}}\eta-v^{p-1}(\iota_v\eta).
\end{equation}
\end{lemma}
\begin{corollary}\label{cor:cartier-pclosed}
Let $\eta$ be a nonzero closed rational one-form on $X$, and let $\F$ be the
foliation defined by $\eta$. Then $\F$ is $p$-closed if and only if
\[
\eta\wedge\Cart_X(\eta)=0.
\]
Equivalently, $\Cart_X(\eta)=a\eta$ for some $a\in k(X)$, possibly zero.
\end{corollary}
\begin{proof}
For every rational vector field $v$ tangent to $\F$, one has $\iota_v\eta=0$.
Formula~\eqref{eq:cartier-contraction} therefore becomes
\[
\bigl(\iota_v\Cart_X(\eta)\bigr)^p=\iota_{v^{[p]}}\eta.
\]
Thus $v^{[p]}$ is tangent to $\F$ for every such $v$ if and only if $\Cart_X(\eta)$
annihilates $\F$ at the generic point. Since $\F$ has corank one, this is equivalent
to $\eta\wedge\Cart_X(\eta)=0$. The conclusion on $X$ follows from
Lemma~\ref{lem:open}.
\end{proof}

\subsection{Smooth pull-back and descent}
We will use the following form of smooth descent.
\begin{lemma}\label{lem:smooth-descent}
Let $q:X\to Y$ be a smooth surjective projective morphism with connected fibres
between smooth varieties. Suppose that $\F\subset T_X$ contains $T_{X/Y}$.
Then there is a foliation $\G\subset T_Y$ with $\F=q^{-1}\G$, and
\begin{equation}\label{eq:smooth-canonical}
K_{\F}\sim K_{X/Y}+q^*K_{\G}.
\end{equation}
If $\F$ has canonical singularities, then $\G$ has canonical singularities.
If $\F$ admits local holomorphic first integrals at every point, then so does $\G$.
\end{lemma}
\begin{proof}
The descent of $\F$ follows from \cite[Lemma~6.7]{AD13}. By construction, on the inverse image of the regular locus of $\G$ there is an exact sequence
\[
0\longrightarrow T_{X/Y}\longrightarrow\F\longrightarrow q^*\G\longrightarrow0.
\]
The complement has codimension at least two. Taking determinants proves
\eqref{eq:smooth-canonical}.

For singularities, take a projective birational morphism
$\mu:Y'\to Y$ with $Y'$ smooth and form $X'=X\times_Y Y'$. Let $\pi:X'\to X$
and $q':X'\to Y'$ be the projections. Then $X'$ is smooth and the relative
canonical bundle commutes with this base change. Writing
$\G'=\mu^{-1}\G$ and $\F'=\pi^{-1}\F$, we obtain
\begin{equation}\label{eq:discrepancy-descent}
K_{\F'}-\pi^*K_{\F}=q'^*(K_{\G'}-\mu^*K_{\G}).
\end{equation}
Since the base change is smooth, each $\mu$-exceptional prime divisor
pulls back with multiplicity one to $\pi$-exceptional prime divisors
with the same discrepancy. As $\mu$ is arbitrary, canonicality
descends from $\F$ to $\G$.

For the local-first-integral assertion, choose $x\in X$ over
$y\in Y$. In a local product neighbourhood for the smooth map $q$,
a local holomorphic first integral $h$ of $\F=q^{-1}\G$ is constant
along the fibres. Thus $h=g\circ q$ for a nonconstant holomorphic
function $g$ on the base, with $dg|_{\G}=0$.
Hence $g$ is a local holomorphic first integral of $\G$ at $y$.
\end{proof}

\subsection{Foliations with numerically trivial canonical class}
We recall the known case of Conjecture~\ref{conj:esbt} for corank one
foliations with canonical singularities and $K_{\F}\equiv0$.
\begin{theorem}[\cite{Dru21}; \cite{DO22}]\label{thm:numtriv}
Let $Z$ be a normal projective klt variety and let $\G\subset T_Z$ be a
positive-rank corank one foliation with canonical singularities. If $K_{\G}\equiv0$
and $\G$ is $p$-closed for almost all primes, then $\G$ is algebraically integrable.
\end{theorem}

\begin{corollary}\label{cor:pure-numtriv}
Let $Z$ be projective and klt. A positive-rank corank one purely transcendental
foliation $\G$ with $K_{\G}$ $\QQ$-Cartier and $K_{\G}\equiv0$ cannot be
$p$-closed for almost all primes.
\end{corollary}
\begin{proof}
By \cite[Proposition~4.22]{Dru21}, a foliation with numerically trivial canonical
class is canonical if and only if it is not uniruled. Pure transcendence excludes a
covering family of rational curves tangent to $\G$, so $\G$ is canonical.
Theorem~\ref{thm:numtriv} implies that $\G$ is algebraically integrable, contradicting
positive rank and pure transcendence.
\end{proof}
\section{Closed rational one-forms and transversely affine foliations}
We first prove algebraic integrability for foliations defined by
closed rational one-forms in Theorem~\ref{thm:closed}.
We then show Corollary~\ref{thm:conormal} for foliations satisfying
$\kappa(N_{\F}^*(D))\geq0$, where $D$ is a reduced simple normal
crossing divisor whose components are $\F$-invariant.
Finally, we prove Theorem~\ref{main:forms} for transversely affine
foliations.

\subsection{Foliations defined by closed rational 1-forms}
\begin{lemma}\label{lem:poles}
Let $X$ be a smooth variety over a perfect field $k$ of characteristic $p>0$,
let $D=\sum_i D_i$ be a reduced simple normal crossing divisor, and let
$P=\sum_i m_iD_i$ with $1\leq m_i<p$. If $\omega$ is a closed rational one-form
with poles bounded by $P$, then
\[
\Cart_X(\omega)\in H^0(X,\Omega_X^1(\log D)).
\]
At a component of $D$ where $\omega$ is logarithmic, the Cartier operator takes its residue
to the $p$-th root of that residue.
\end{lemma}
\begin{proof}
The assertion is local at the generic points of prime divisors, since
$\Omega_X^1(\log D)$ is locally free. Away from $D$, Cartier sends closed regular
forms to regular forms. Fix the generic point of a component $D_i$, choose a
uniformizer $t$, and put $L=k(D_i)$. In the completed local field $L((t))$, use
a separating transcendence basis $u_2,\ldots,u_n$ of $L/k$ and write
\[
\omega=a(t)\,dt+\sum_{j=2}^n b_j(t)\,du_j.
\]
The Laurent coefficients have no terms of order less than $-m_i$.
Closedness implies
$\partial_t b_j=\partial_{u_j}a$. Since $m_i<p$, the coefficient of
$t^{-m_i}$ in each $b_j$ vanishes: its derivative would have order $-m_i-1$,
which does not occur in $\partial_{u_j}a$.

For $2\leq r\leq m_i$, a term $a_{-r}t^{-r}dt$ can be removed by subtracting
\[
d\left(\frac{a_{-r}}{1-r}t^{1-r}\right).
\]
The denominators $1-r$ are invertible, and the new tangential terms have smaller
pole order. Repeating the procedure and using the same closedness identity shows
\[
\omega=dh+c\,\frac{dt}{t}+\beta,
\]
where $\beta$ is regular in the completed local ring. Closedness implies
$dc=0$ in $\Omega^1_{L/k}$, hence $c\in L^p$, and $\beta$ is closed.
By Proposition~\ref{prop:cartier} and
Corollary~\ref{cor:cartier-rules}, we obtain
\[
\Cart_X(\omega)=c^{1/p}\frac{dt}{t}+\Cart_X(\beta).
\]
The second summand is regular. This proves the logarithmic pole bound at $D_i$
and the residue assertion. Checking all prime divisors proves the lemma.
\end{proof}

We recall the following two results. A flat algebraic connection is
called \emph{isotrivial} if it becomes trivial after a finite etale cover.

\begin{theorem}[\cite{And04}, Corollaries~4.3.6 and~7.1.3]\label{thm:rank-one}
Let $K$ be a field finitely generated over $\QQ$, let $U$ be a smooth
geometrically connected $K$-variety, and let $(L,\nabla)$ be a rank-one
flat algebraic connection on $U$. If its reductions have vanishing
$p$-curvature for almost all primes, then $(L,\nabla)$ is isotrivial.
\end{theorem}

\begin{theorem}[\cite{And04}, Corollaries~4.3.4 and~7.1.3]\label{thm:andre-extension}
Let $K$ be a field finitely generated over $\QQ$, let $U$ be a smooth
geometrically connected $K$-variety, and let
\[
0\longrightarrow(E',\nabla')
\longrightarrow(E,\nabla)
\longrightarrow(E'',\nabla'')
\longrightarrow0
\]
be an exact sequence of flat algebraic connections on $U$.
Suppose that the two outer terms are isotrivial.
If the reductions of $(E,\nabla)$ have vanishing $p$-curvature for
almost all primes, then $(E,\nabla)$ is isotrivial.
\end{theorem}

\begin{lemma}\label{lem:wedge}
Let $X$ be a variety with function field $K=\CC(X)$,
and let $\omega$ be a nonzero closed rational one-form on $X$.
Assume that $\ker\omega$ has no nonconstant rational first integral.
Let $W$ be a finite-dimensional complex vector space of closed
rational one-forms on $X$ containing $\omega$. Then
\begin{equation}\label{eq:wedge-kernel}
\ker\bigl(W\longrightarrow\Omega^2_{K/\CC},
          \ \theta\longmapsto\omega\wedge\theta\bigr)
=\CC\omega.
\end{equation}

Moreover, choose a basis $\omega,\theta_2,\ldots,\theta_r$ of $W$
and spread out $X$ and these forms over
$S=\Spec R$, where $R\subset\CC$ is an integral finitely generated
$\ZZ$-algebra.
After shrinking $S$, the reduced forms
$\omega_s,\theta_{2,s},\ldots,\theta_{r,s}$ are linearly independent
over $k(s)$.
Writing
\[
K_s=k(s)(X_s),
\qquad
W_s=\operatorname{Span}_{k(s)}
\{\omega_s,\theta_{2,s},\ldots,\theta_{r,s}\},
\]
we have
\[
\ker\bigl(W_s\longrightarrow\Omega^2_{K_s/k(s)},
          \ \theta\longmapsto\omega_s\wedge\theta\bigr)
=k(s)\omega_s
\]
for every closed point $s\in S$.
\end{lemma}
\begin{proof}
Let $\theta\in W$ satisfy $\omega\wedge\theta=0$.
Since $\omega\neq0$, we can write $\theta=b\omega$ for some $b\in K$.
Both forms are closed, so
\[
0=d\theta=db\wedge\omega.
\]
If $b$ were nonconstant, it would be a rational first integral of
$\ker\omega$. Thus $b\in\CC$, proving \eqref{eq:wedge-kernel}.

We now explain why this kernel remains one-dimensional after reduction.
Write $r=\dim_{\CC}W$ and choose a basis
\(
\omega,\theta_2,\ldots,\theta_r
\)
of $W$. Choose a smooth projective model $Y$ of $K$ and an effective
divisor $P$ bounding the poles of these forms. Then
\[
W\subset H^0(Y,\Omega_Y^1(P)),
\qquad
\omega\wedge W\subset H^0(Y,\Omega_Y^2(2P)).
\]
Extend the chosen basis of $W$ to a basis of the first space, and
choose a basis of the second space.

Spread out $Y$, $P$, and these sections over an integral finitely
generated $\ZZ$-subalgebra $R\subset\CC$.
After shrinking $\Spec R$, cohomology and base change ensure that
the chosen bases reduce to bases of the corresponding spaces of
sections on every closed fibre.
In particular,
\[
W_s=\operatorname{Span}_{k(s)}
\{\omega_s,\theta_{2,s},\ldots,\theta_{r,s}\}
\]
has dimension $r$, and $\omega_s\neq0$.

With respect to these bases, the wedge map on $W$ is represented
by a matrix $M$ with entries in $R$.
Its first column is zero because $\omega\wedge\omega=0$, and
\eqref{eq:wedge-kernel} shows that its rank is $r-1$.
If $r=1$, the conclusion is immediate.
Otherwise, choose a nonzero $(r-1)\times(r-1)$ minor $\delta$ of $M$
and replace $R$ by $R[\delta^{-1}]$.
On every remaining closed fibre, this minor stays nonzero, so
$\operatorname{rank}M_s\geq r-1$.
The first column of $M_s$ is still zero, so
$\operatorname{rank}M_s\leq r-1$.
Consequently,
\[
\dim_{k(s)}\ker M_s=1.
\]
Since $\omega_s$ is nonzero and belongs to this kernel, we obtain
\[
\{\theta\in W_s:\omega_s\wedge\theta=0\}=k(s)\omega_s,
\]
as required.
\end{proof}

\begin{theorem}\label{thm:closed}
Let $X$ be a normal complex variety of dimension at least two, and let $\F$ be a
corank one foliation defined by a nonzero closed rational one-form. If $\F$ is
$p$-closed for almost all primes, then it has a nonconstant rational first integral.
In particular, $\F$ is algebraically integrable.
\end{theorem}
\begin{proof}
Let $\omega$ be a nonzero closed rational one-form defining $\F$. The existence of a rational first integral is birationally invariant,
as is the property of being $p$-closed for almost all primes by
Lemma~\ref{lem:open}.
By compactifying a dense affine open subset and resolving, we may
therefore assume that $X$ is smooth and projective and that the support
$D$ of the polar divisor $P=\sum_i m_iD_i$ of $\omega$ is simple normal
crossing. Suppose, for a contradiction, that $\F$ is not algebraically
integrable, or equivalently that it has no nonconstant rational first
integral.

Set
\begin{equation}\label{eq:log-space}
V=H^0(X,\Omega_X^1(\log D)),
\qquad
W=V+\CC\omega.
\end{equation}
Every element of $V$ is closed by \cite[Corollary~5.4]{Per24}, so
Lemma~\ref{lem:wedge} applies to $W$.
Choose a model over $S=\Spec R$ including the pair $(X,D)$,
the pole bound $P$, the form $\omega$, and bases of $V$ and $W$.
After shrinking $S$, we may assume that the model of $X$ is smooth
and projective over $S$, with relative simple normal crossing boundary,
and that cohomology and base change identify
\[
V_s=H^0(X_s,\Omega_{X_s}^1(\log D_s)).
\]
We may also assume that
\begin{equation}\label{eq:log-kernel}
\{\theta\in W_s:\omega_s\wedge\theta=0\}=k(s)\omega_s
\end{equation}
on every closed fibre, that the reduced foliations are $p(s)$-closed,
and that $p(s)>\max_i m_i$ whenever $D\neq0$.
Put $U=X\setminus D$.

\noindent \textbf{Claim.} $\omega\in V$.

\begin{proof}

Suppose first that $\omega\notin V$.
After shrinking the base further, we may assume that
$\omega_s\notin V_s$ on every closed fibre.
Lemma~\ref{lem:poles} shows $\Cart(\omega_s)\in V_s$, while
Corollary~\ref{cor:cartier-pclosed} implies
$\omega_s\wedge\Cart(\omega_s)=0$.
Equation~\eqref{eq:log-kernel} therefore implies
\[
\Cart(\omega_s)\in V_s\cap k(s)\omega_s=\{0\}.
\]

Since $\omega|_U$ is a closed regular one-form, the trivial rank-two
vector bundle
\[
E=\OO_Ue_1\oplus\OO_Ue_2
\]
admits the flat algebraic connection defined by
\[
\nabla e_1=0,
\qquad
\nabla e_2=-e_1\otimes\omega.
\]
It fits into an exact sequence of flat connections
\[
0\longrightarrow(\OO_U,d)
\longrightarrow(E,\nabla)
\longrightarrow(\OO_U,d)
\longrightarrow0,
\]
where the subbundle is generated by $e_1$ and the quotient is generated
by the image of $e_2$.

Fix a closed point $s$ and write $p=p(s)$.
For a local vector field $v$ on $U_s$, the $p$-curvature
\[
\psi_{\nabla_s}(v)
=\nabla_{s,v}^{\,p}-\nabla_{s,v^{[p]}}
\]
satisfies $\psi_{\nabla_s}(v)(e_1)=0$ and
\[
\begin{aligned}
\psi_{\nabla_s}(v)(e_2)
&=\bigl(\omega_s(v^{[p]})
       -v^{p-1}(\omega_s(v))\bigr)e_1\\
&=\bigl(\iota_v\Cart(\omega_s)\bigr)^p e_1
=0,
\end{aligned}
\]
by \eqref{eq:cartier-contraction}.
Thus $(E,\nabla)$ has vanishing $p$-curvature for almost all primes.

Both outer terms are trivial, so
Theorem~\ref{thm:andre-extension} implies a connected finite \'etale cover
$\pi:U'\to U$ trivializing $(E,\nabla)$.
Let $\rho:\pi^*E\to\OO_{U'}$ be the quotient map, and choose an
algebraic horizontal frame $s_1,s_2$ of $\pi^*E$.
Since $\rho$ respects the connections, $d(\rho(s_i))=0$,
so $\rho(s_1)$ and $\rho(s_2)$ are constants.
Surjectivity of $\rho$ implies that at least one is nonzero.
After relabelling and rescaling, we may therefore assume
$\rho(s_1)=1$.

In the pulled-back frame $e_1,e_2$, this means that
$s_1=e_2+he_1$ for some $h\in H^0(U',\OO_{U'})$.
Its horizontality implies
\[
0=(\pi^*\nabla)(s_1)
=e_1\otimes(dh-\pi^*\omega),
\]
hence $dh=\pi^*\omega$.
Writing $e=[\CC(U'):\CC(X)]$, the rational function
\[
f=\frac{1}{e}\Tr_{\CC(U')/\CC(X)}(h)
\]
satisfies $df=\omega$.
Since $\omega\neq0$, this is a nonconstant rational first integral,
a contradiction.
\end{proof}
So we have $\omega\in V$, and 
$\Cart(\omega_s)\in V_s$.
Together with Corollary~\ref{cor:cartier-pclosed} and
\eqref{eq:log-kernel}, this yields
\begin{equation}\label{eq:constant-cartier}
\Cart(\omega_s)=c_s\omega_s,
\qquad c_s\in k(s).
\end{equation}
\noindent \textbf{Claim.} $\omega$ is holomorphic.
\begin{proof}

Suppose that the residue $\lambda$ of $\omega$ along some component
of $D$ is nonzero. Set
\(
\beta=\frac{\omega}{\lambda},
\)
and include $\lambda$ and $\lambda^{-1}$ in the ring of definition.
The residue of $\beta_s$ along the corresponding component is then $1$.
By Frobenius semilinearity, \eqref{eq:constant-cartier} implies
\[
\Cart(\beta_s)=c'_s\beta_s
\qquad\text{for some }c'_s\in k(s).
\]
Taking residues and applying Lemma~\ref{lem:poles} implies $c'_s=1$.
Thus
\[
\Cart(\beta_s)=\beta_s.
\]

On the trivial line bundle $L=\OO_U$, consider the flat algebraic
connection
\[
\nabla=d-\beta.
\]
For a local vector field $v$ on $U_s$, with $p=p(s)$, its
$p$-curvature satisfies
\[
\begin{aligned}
\psi_{\nabla_s}(v)(1)
&=\beta_s(v^{[p]})
  -v^{p-1}(\beta_s(v))-\beta_s(v)^p\\
&=\bigl(\iota_v\Cart(\beta_s)\bigr)^p-\beta_s(v)^p
=0,
\end{aligned}
\]
again by \eqref{eq:cartier-contraction}.
Theorem~\ref{thm:rank-one} therefore implies the existence of a connected finite
\'etale cover $\pi:U'\to U$ on which $(L,\nabla)$ has a
nowhere-vanishing algebraic horizontal section $h$.
The horizontal equation is
\[
\frac{dh}{h}=\pi^*\beta.
\]
Taking the norm, we obtain
\[
f=\operatorname{Nm}_{\CC(U')/\CC(X)}(h)\in\CC(X)^\times,
\qquad
\frac{df}{f}=e\beta,
\]
where $e=[\CC(U'):\CC(X)].$
Since $\beta\neq0$, the function $f$ is nonconstant and is a rational
first integral of $\F$, again a contradiction. Thus all residues of $\omega$ along $D$ vanish.
Since $\omega\in H^0(X,\Omega_X^1(\log D))$, it follows that
$\omega$ has no poles and is therefore holomorphic.
\end{proof}
Let $a:X\to A=\Alb(X)$ be an Albanese morphism.
There is a unique invariant one-form $\alpha$ on $A$ such that
$a^*\alpha=\omega$.

If $\dim A=1$, then $a$ is a nonconstant morphism to a curve inducing
$\F$. Composing it with a nonconstant rational function yields a
rational first integral, a contradiction.
We may therefore assume that $\dim A>1$. After spreading out $X$, $A$, $a$, and $\alpha$ over a common base
and shrinking it, we may assume that $A$ extends to an abelian
scheme and that
\(
a_s^*:H^0(A_s,\Omega_{A_s}^1)
\rightarrow H^0(X_s,\Omega_{X_s}^1)
\)
is an isomorphism on every closed fibre.
Functoriality of Cartier and \eqref{eq:constant-cartier} implies
\[
a_s^*\bigl(\Cart(\alpha_s)-c_s\alpha_s\bigr)=0.
\]
Injectivity of $a_s^*$ therefore implies
\[
\Cart(\alpha_s)=c_s\alpha_s.
\]
By Corollary~\ref{cor:cartier-pclosed}, the invariant foliation
$\HH=\ker\alpha$ on $A$ is $p$-closed for almost all primes.
By \cite[Proposition~10.7]{Dru21}, $\HH$ is algebraically integrable.
Since $\HH$ is translation invariant, it is induced by a quotient
morphism $q:A\to E$ to an elliptic curve, and $\alpha=q^*\beta$
for a nonzero holomorphic one-form $\beta$ on $E$.
The equality
\[
(q\circ a)^*\beta=a^*\alpha=\omega\neq0
\]
shows that $q\circ a$ is nonconstant and induces $\F$.
Composing $q\circ a$ with a nonconstant rational function on $E$
defines a rational first integral of $\F$, a contradiction.
This completes the proof.
\end{proof}

\begin{corollary}\label{thm:conormal}
Let $X$ be a smooth complex projective variety and let $\F$ be a positive-rank
corank one foliation. Let $D$ be a reduced simple normal crossing divisor whose
irreducible components are $\F$-invariant.
Assume that $\F$ is $p$-closed for almost all primes and that
$\kappa(N_{\F}^*(D))\geq0$. Then $\F$ is algebraically integrable.
\end{corollary}
\begin{proof}
Set $L=N_{\F}^*(D)$. Invariance of every component of $D$ implies an inclusion
\[
L\longrightarrow\Omega_X^1(\log D).
\]
Indeed, if $D$ has local equation $t_1\cdots t_r=0$ and $\omega$ is a local
generator of the conormal image in $\Omega_X^1$, then invariance says that
$t_i$ divides $\omega\wedge dt_i$. Consequently
$\omega/(t_1\cdots t_r)$ is logarithmic along $D$.

Choose an integer $m>0$ and a nonzero section
$s\in H^0(X,L^{\otimes m})$.
Take the cyclic cover obtained by extracting an $m$-th root of $s$. Choose an irreducible component dominating $X$,
normalize it, and take a projective resolution on which the reduced inverse
image $E$ of $D$ is simple normal crossing. We obtain a smooth projective
variety $Y$ and a generically finite morphism $\pi:Y\to X$.
The tautological root defines a nonzero section
$\tau\in H^0(Y,\pi^*L)$ satisfying $\tau^{\otimes m}=\pi^*s$. Since $\Supp(\pi^*D)\subseteq E$, pull-back of differential forms
induces a natural morphism
\[
\pi^*\Omega_X^1(\log D)\longrightarrow\Omega_Y^1(\log E).
\]

Let $\eta\in H^0(Y,\Omega_Y^1(\log E))$ be the image of $\tau$ under
the composite morphism
\[
\pi^*L\longrightarrow\pi^*\Omega_X^1(\log D)
\longrightarrow\Omega_Y^1(\log E).
\]
Since $\pi$ is generically separable, pull-back is injective on rational
one-forms. At the generic point of $Y$, the section $\tau$ generates
$\pi^*L$, whose image is the conormal line of the pulled-back foliation.
Thus $\eta$ is nonzero and defines $\pi^{-1}\F$, that is,
\[
(\ker\eta)^{\sat}=\pi^{-1}\F.
\]

Since $Y$ is smooth and projective and $E$ has simple normal crossing
support, the global logarithmic one-form $\eta$ is closed by
\cite[Corollary~5.4]{Per24}.
Moreover, Lemma~\ref{lem:finite} implies that $\pi^{-1}\F$ is
$p$-closed for almost all primes.
Theorem~\ref{thm:closed}, applied to the foliation defined by $\eta$,
therefore proves that $\pi^{-1}\F$ is algebraically integrable.
Algebraic integrability of $\F$ follows by descent along $\pi$,
using Lemma~\ref{lem:finite}.
\end{proof}

\begin{remark}
By \cite[Corollary~9.2]{Tou16},
$\kappa(N_{\F}^*(D))\leq1$, and equality implies algebraic
integrability.
The corollary proves algebraic integrability
when $\kappa(N_{\F}^*(D))=0$ when $\F$ is
$p$-closed for almost all primes.
\end{remark}
\subsection{Transversely affine foliations}
\begin{definition}\label{def:affine}
Let $X$ be a normal complex variety and let $\F$ be a corank one
foliation on $X$. We say that $\F$ is \emph{transversely affine}
if there exist rational one-forms $\omega\neq0$ and $\eta$ such that
\[
\F=\ker\omega,\qquad
d\omega=\omega\wedge\eta,\qquad d\eta=0.
\]

When $X$ is smooth and projective, this is equivalent to the
existence of a meromorphic flat connection $\nabla$ on $N_{\F}$
satisfying $d_{\nabla}\omega_{\F}=0$, where $\omega_{\F}$ is the
defining twisted one-form of Remark~\ref{omegaF} and $d_{\nabla}$
is the exterior covariant derivative.
\end{definition}
\begin{definition}
    A Riccati foliation on a ruled surface
$q:S\to B$ is a rank-one foliation transverse to the general fibre of $q$;
over a nonempty open subset of $B$, it is the horizontal foliation of a
flat $\PP^1$-connection.
\end{definition}
We recall the following structure theorem in the form in which it will be used.
\begin{theorem}[{\cite[Theorem~2.2]{CLNP22}}]\label{thm:affine-structure}
Let $X$ be a smooth complex projective variety and let $\F$ be a
transversely affine corank one foliation. Then at least one of the
following holds:
\begin{enumerate}
\item there is a generically finite cover $\pi:Y\to X$ such that
$\pi^{-1}\F$ is defined by a closed rational one-form;
\item $\F$ is the rational pull-back of a transversely affine Riccati
foliation on a ruled projective surface.
\end{enumerate}
\end{theorem}

\begin{lemma}\label{lem:affine-riccati}
Let $\mathcal R$ be a transversely affine Riccati foliation on a smooth ruled
projective surface $q:S\to B$. If $\mathcal R$ is $p$-closed for almost all
primes, then it is algebraically integrable.
\end{lemma}
\begin{proof}
Let $\nabla$ be a meromorphic flat connection on $N_{\mathcal R}$
defining the transverse affine structure, and let $D$ be the reduced
divisor supported on the poles of $\nabla$.
Suppose that $D$ is vertical over $B$. Then $q(D)$ is a finite
subset of $B$. Choose $b\in B\setminus q(D)$ such that the fibre
$F=q^{-1}(b)\simeq\PP^1$ is everywhere transverse to $\mathcal R$.
The natural morphism $T_F\to N_{\mathcal R}|_F$ is therefore an
isomorphism, giving
\[
N_{\mathcal R}|_F\simeq T_F\simeq\OO_{\PP^1}(2).
\]
Since $F\cap D=\varnothing$, the restriction of $\nabla$ to $F$
is a holomorphic connection on $\OO_{\PP^1}(2)$.
This is impossible, since a line bundle on a smooth projective
curve admitting a holomorphic connection has degree zero.
Thus $D$ has an irreducible component dominating $B$.
 By \cite[Proposition~2.1]{CP14},
every component of $D$ is $\mathcal R$-invariant. We obtain an invariant
multisection of $q$.

Let $\nu:B'\to C$ be the normalization of the invariant multisection
$C\subset S$. The composite $\tau=q\circ\nu:B'\to B$ is finite.
The base change
\[
q':S':=S\times_B B'\longrightarrow B'
\]
has a natural section
\[
\sigma:B'\longrightarrow S',
\qquad b'\longmapsto(\nu(b'),b').
\]
By the Riccati property, after removing finitely many points from
$B'$, we obtain an open curve $B^\circ$ on which $\tau$ is \'etale,
$q'$ is a $\PP^1$-bundle, and the pulled-back foliation is regular
and everywhere transverse to the fibres.
The section $\sigma(B^\circ)$ is invariant because $C$ is invariant,
and the pulled-back foliation defines a flat horizontal connection
on this $\PP^1$-bundle.

Denote the restricted $\PP^1$-bundle by $q^\circ:P\to B^\circ$
and its invariant section by $\Sigma\subset P$. Set
\[
E=q^\circ_*\OO_P(\Sigma).
\]
Since $\OO_P(\Sigma)$ restricts to $\OO_{\PP^1}(1)$ on each
fibre, cohomology and base change show that $E$ is locally free
of rank two. The evaluation map
\[
(q^\circ)^*E\longrightarrow\OO_P(\Sigma)
\]
is surjective because $\OO_{\PP^1}(1)$ is globally generated.
It therefore induces a morphism $\varphi:P\to\PP(E)$ over
$B^\circ$, where $\PP(E)$ parametrizes one-dimensional quotients.
On each fibre, $\varphi$ is the isomorphism defined by the
complete linear system $|\OO_{\PP^1}(1)|$.
Hence $\varphi$ is an isomorphism.

Let $\mathcal R^\circ$ denote the restriction to $P$ of the
pull-back of $\mathcal R$ to $S'$. This foliation is regular
and transverse to the fibres of $q^\circ$.
We use $\mathcal R^\circ$ to construct an algebraic connection
on $E$.

Let $V\subset B^\circ$ be open. By the definition of
$E=q^\circ_*\OO_P(\Sigma)$, a section $f\in H^0(V,E)$ is a section
of $\OO_P(\Sigma)$ on $(q^\circ)^{-1}(V)$, and can therefore be
viewed as a rational function on this inverse image with poles
bounded by $\Sigma$.

For a local vector field $v$ on $V$, let $\widetilde v$ be the
unique vector field on $(q^\circ)^{-1}(V)$ tangent to
$\mathcal R^\circ$ and satisfying $dq^\circ(\widetilde v)=v$.
We call $\widetilde v$ the horizontal lift of $v$ with respect
to $\mathcal R^\circ$. Since $\Sigma$ is invariant, if a local rational function $f$
has poles bounded by $\Sigma$, then $\widetilde v(f)$ also has
poles bounded by $\Sigma$. Thus differentiation by $\widetilde v$
sends local sections of $\OO_P(\Sigma)$ to local sections of
the same sheaf. We may therefore define
\[
\nabla_{E,v}(f)=\widetilde v(f).
\]

Thus $\nabla_E$ is an algebraic
connection on $E$. It is flat because $B^\circ$ is a curve.

By construction, the isomorphism $\varphi:P\to\PP(E)$ identifies
$\mathcal R^\circ$ with the horizontal foliation induced by
$\nabla_E$ on $\PP(E)$.
The constant function $1$ defines a nowhere vanishing horizontal
section of $E$. Thus there is an exact sequence of flat connections
\[
0\longrightarrow(\OO_{B^\circ},d)
\longrightarrow(E,\nabla_E)
\longrightarrow(L,\nabla_L)
\longrightarrow0,
\qquad \rank L=1.
\]

Spread out these data and shrink the arithmetic base so that the
construction above commutes with reduction.
By Lemma~\ref{lem:finite}, we may also assume that the reduced
horizontal foliations are $p$-closed on every closed fibre.

Fix a closed point $s$ of the arithmetic base and put
$p=\operatorname{char}k(s)$.
For a local vector field $v$ on $B^\circ_s$, let $\widetilde v$
be its horizontal lift to $P_s$.
Both $\widetilde v^{[p]}$ and $\widetilde{v^{[p]}}$ project
to $v^{[p]}$, so their difference is vertical.
It is also horizontal by $p$-closedness, and therefore vanishes.
Consequently, for every local section $f$ of $E_s$,
\[
\begin{aligned}
\psi_{\nabla_{E,s}}(v)(f)
&=\nabla_{E,s,v}^{\,p}(f)
  -\nabla_{E,s,v^{[p]}}(f)\\
&=\bigl(\widetilde v^{[p]}
        -\widetilde{v^{[p]}}\bigr)(f)
=0.
\end{aligned}
\]
Thus $(E,\nabla_E)$ has zero $p$-curvature for almost all primes.

The rank-one quotient $(L,\nabla_L)$ has zero $p$-curvature for
almost all primes, since it is a quotient of $(E,\nabla_E)$.
The rank-one Grothendieck--Katz theorem
\cite[Corollaries~4.3.6 and~7.1.3]{And04}
therefore implies that $(L,\nabla_L)$ is isotrivial.
Since the subconnection $(\OO_{B^\circ},d)$ is trivial,
Andr\'e's extension theorem
\cite[Corollaries~4.3.4 and~7.1.3]{And04}
implies that $(E,\nabla_E)$ is isotrivial.

After a finite \'etale base change, its projectivization is the
trivial $\PP^1$-bundle with the product connection.
Its horizontal leaves are therefore algebraic.
Algebraic integrability descends by Lemma~\ref{lem:finite},
first through this cover and then through the original finite
base change.
\end{proof}

\begin{theorem}\label{thm:affine}
Let $X$ be a normal complex variety of dimension at least two,
and let $\F$ be a transversely affine corank one foliation.
If $\F$ is $p$-closed for almost all primes, then $\F$ is
algebraically integrable.
\end{theorem}
\begin{proof}
By Lemma~\ref{lem:open}, we may replace $X$ by a smooth projective
model, since both the transverse affine structure and algebraic
integrability are birationally invariant.

We apply Theorem~\ref{thm:affine-structure}.
In the first case, there is a generically finite cover
$\pi:Y\to X$ such that $\pi^{-1}\F$ is defined by a closed rational
one-form. After resolving $Y$, we may assume that $Y$ is smooth.
By Lemma~\ref{lem:finite}, $\pi^{-1}\F$ is $p$-closed for almost
all primes, so Theorem~\ref{thm:closed} implies that it is
algebraically integrable. Algebraic integrability then descends
to $\F$ by Lemma~\ref{lem:finite}.

In the second case, there is a rational map $f:X\dashrightarrow S$
to a ruled surface and a transversely affine Riccati foliation
$\mathcal R$ on $S$ such that $\F=f^{-1}\mathcal R$.
If the image of $f$ is a curve, then $\F$ is induced by the
rational map from $X$ to the normalization of that curve,
and hence is algebraically integrable.
We may therefore assume that $f$ is dominant.
By Lemma~\ref{lem:dominant}, $\mathcal R$ is $p$-closed for almost
all primes. Lemma~\ref{lem:affine-riccati} implies that
$\mathcal R$ is algebraically integrable, and therefore so is
$\F=f^{-1}\mathcal R$.
\end{proof}

This proves Theorem~\ref{main:forms}.
In particular, on a smooth projective variety the theorem also
applies to a defining meromorphic one-form with values in a
line bundle which is closed for a specified meromorphic flat
connection.

\section{Foliations with nef anti-canonical divisor}
We recall two results on foliations with nef anti-canonical class
that will be used in the proof of Theorem~B.

\begin{proposition}[{\cite[Proposition~6.1(2)]{Dru17};
\cite{DruErr}}]\label{prop:comparison}
Let $X$ be smooth and projective, let $\F\subset T_X$ be a foliation,
and let $\HH=\F^a$. If $-K_{\F}$ is nef and $\HH$ has a compact leaf,
then $K_{\HH}\equiv K_{\F}$.
\end{proposition}

\begin{theorem}[{\cite[Theorem~1.1 and Corollary~1.2]{Ou25}}]
\label{thm:ou}
Let $\HH$ be a foliation on a smooth projective variety $X$.
Assume that $-K_{\HH}$ is nef and that $\HH$ has a compact leaf.
Then there is a smooth, locally trivial projective fibration
$q:X\to Y$ with rationally connected fibres and a foliation $\HH_0$
on the smooth projective variety $Y$ such that
\[
\HH=q^{-1}\HH_0,\qquad K_{\HH_0}\equiv0.
\]
Moreover, $\HH$ is regular, and there is a regular foliation $\E$
such that
\[
T_X=\HH\oplus\E.
\]
\end{theorem}
\subsection{Local holomorphic first integrals and algebraic parts}
The next lemma is an extension of the proof of \cite[Lemma~6.2]{Dru17}, which assumes regularity or the existence of a compact leaf of the foliation.
\begin{lemma}\label{lem:local-compact}
Let $X$ be smooth and projective, and let $\F$ be a corank one foliation
admitting local holomorphic first integrals at every point. If $r^a(\F)>0$,
then its algebraic part $\HH=\F^a$ is induced by an almost holomorphic map.
In particular, $\HH$ has a compact leaf.
\end{lemma}
\begin{proof}
Let $Z$ be the normalization of the family of closures of general leaves of
$\HH$, and let $U$ be the normalization of its universal cycle.
We have projective morphisms
\[
q:U\longrightarrow Z,\qquad \nu:U\longrightarrow X,
\]
where $\nu$ is birational and $(q,\nu):U\to Z\times X$ is finite; see
\cite[Section~2.2, paragraph~2.13]{Dru17}. 
Let $\G$ be the purely transcendental quotient on $Z$, so that
$\nu^{-1}\F=q^{-1}\G$. If $\G=0$, the same notation is used.\\

\noindent \textbf{Claim.} $q(\Exc(\nu))\neq Z$. 
\begin{proof}
Suppose for contradiction that $q(\operatorname{Exc}(\nu))=Z$. Since $X$ is smooth, the exceptional locus of $\nu$ is a union
of exceptional prime divisors. Suppose an exceptional
prime divisor $E$ dominates $Z$.

By generic smoothness,
we can choose a sufficiently general point $u\in E$ over a very
general point $z\in Z$ such that $U$ and $E$ are smooth at $u$,
$Z$ and $\G$ are regular at $z$, and $q|_E$ is smooth at $u$.
In particular, the differential of $q$ is surjective at $u$,
so $q$ is a submersion near $u$.

Since $u$ is a general point of the exceptional divisor $E$,
the projective fibre of $\nu$ through $u$ is positive-dimensional
at $u$.
Choose a complete integral curve $\Gamma$ in this fibre through $u$.
Finiteness of $(q,\nu)$ implies that $q|_\Gamma$ is finite.
Consequently, $B=q(\Gamma)$ is a curve, and a local branch of
$\Gamma$ through $u$ maps nonconstantly to a local branch of $B$
through $z$.

Put $x=\nu(u)$.
We show that the pull-back of a local first integral of $\F$
at $x$ descends to a nonconstant local first integral of $\G$
at $z$.
Choose a nonconstant holomorphic first integral $h$ of $\F$
on a connected neighbourhood $\Omega$ of $x$.
Since $q$ is a submersion at $u$, we can choose neighbourhoods
$W$ of $u$ and $V$ of $z$ such that $\nu(W)\subset\Omega$ and
$q|_W:W\to V$ is a holomorphic submersion with connected fibres.
We also assume that $\G$ is regular on $V$.

The function $H=h\circ\nu$ is a holomorphic first integral of
$\nu^{-1}\F=q^{-1}\G$ on $W$.
Since the fibres of $q|_W$ are connected and tangent to
$q^{-1}\G$, the function $H$ is constant on each fibre.
It therefore descends to a function $g$ on $V$, with
\begin{equation}\label{eq:local-factor}
H=h\circ\nu=g\circ q
\qquad\text{on }W.
\end{equation}
To see that $g$ is holomorphic, choose a local holomorphic section
$\sigma:V_0\to W$ of $q|_W$.
Since $H$ is constant on each fibre, we have
$g|_{V_0}=H\circ\sigma$, which is holomorphic.

The germ of $g$ at $z$ is nonconstant. Otherwise $h\circ\nu$
would be constant on a neighbourhood of $u$. This neighbourhood
meets the isomorphism locus of $\nu$, so $h$ would be constant
on a nonempty open subset of $\Omega$, and hence on $\Omega$,
a contradiction. Then $g$ is a local holomorphic first integral of $\G$.
Indeed, since $q$ is a submersion, every tangent vector
$\xi\in\G_{q(w)}$ lifts to a vector
$\widetilde\xi\in(q^{-1}\G)_w$.
Thus \eqref{eq:local-factor} implies
\[
dg_{q(w)}(\xi)=dH_w(\widetilde\xi)=0.
\]
Hence $dg|_{\G}=0$.

By \eqref{eq:local-factor}, the equality $H=h(x)$ on
$\Gamma\cap W$ implies that $g=h(x)$ on the corresponding
local branch of $B$ through $z$.
Since $\G$ is regular of corank one and $g$ is a nonconstant
local first integral, the reduced germ of the level set
$\{g=g(z)\}$ at $z$ is the local leaf of $\G$ through $z$.
Hence this branch of $B$ is contained in that local leaf.

So $B$ is tangent to $\G$ on a nonempty analytic
open subset.
Since $B$ is irreducible and tangency is an algebraic condition
on the common regular locus, $B$ is tangent to $\G$ at its
generic point.
This contradicts pure transcendence, as $B$ passes through
the very general point $z$.

\end{proof}

After removing $q(\Exc(\nu))$ from $Z$, the map $\nu$ identifies the universal
family with an open subset of $X$, and $q$ defines a proper morphism there.
Shrinking the base further makes this morphism smooth with connected fibres.
Its fibres are smooth and projective, and its relative tangent bundle is
$\HH$. They are compact leaves of $\HH$ as required.
\end{proof}
\subsection{Proof of Theorem \ref{main:nef}}
\begin{proof}[Proof of Theorem B]
Set $\HH=\F^a$.
Suppose first that $\rank\HH=0$.
By \cite[Corollary~4.13]{CP19}, the canonical class $K_{\F}$ is
pseudo-effective: otherwise $\F$ would contain a positive-rank
algebraically integrable subfoliation.
Since $-K_{\F}$ is nef, it follows that $K_{\F}\equiv0$.
Moreover, $\F$ is not uniruled, because its algebraic part has rank zero.
Thus \cite[Corollary~3.8]{LPT18} implies that $\F$ has canonical
singularities.
Theorem~\ref{thm:numtriv} then makes $\F$ algebraically integrable,
contradicting $\rank\HH=0$.
Hence $\HH$ has positive rank.

The foliation $\HH$ has a compact leaf by
Lemma~\ref{lem:canonical-compact} in the canonical case, and by
Lemma~\ref{lem:local-compact} in the local-first-integral case.
Proposition~\ref{prop:comparison} implies
\begin{equation}\label{eq:comparison}
K_{\HH}\equiv K_{\F}.
\end{equation}
In particular, $-K_{\HH}$ is nef.
Theorem~\ref{thm:ou} implies that $\HH$ is regular and there exists a smooth
projective fibration $q:X\to Y$ with rationally connected fibres
and a foliation $\HH_0$ on $Y$ such that
\[
\HH=q^{-1}\HH_0,
\qquad
K_{\HH_0}\equiv0.
\]
Since $q$ is smooth and $\HH$ is regular, $\HH_0$ is also regular.
Moreover, $T_{X/Y}\subset\HH\subset\F$, so
Lemma~\ref{lem:smooth-descent} implies
$\F=q^{-1}\G$ for a foliation $\G$ on $Y$.
Its corank is one, with the possibility $\G=0$ if $Y$ is a curve.

For the two smooth pull-backs we have
\begin{equation}\label{eq:two-pullbacks}
\begin{aligned}
K_{\F}&\sim K_{X/Y}+q^*K_{\G},\\
K_{\HH}&\sim K_{X/Y}+q^*K_{\HH_0}.
\end{aligned}
\end{equation}
So by \eqref{eq:comparison}, \eqref{eq:two-pullbacks} and $K_{\HH_0}\equiv0$, we have $q^*K_{\G}\equiv0$. Hence $K_{\G}\equiv0$.

If $\G=0$, then $\F=T_{X/Y}$ is algebraically integrable.
Suppose now that $\G\neq0$.
We claim that its algebraic part is $\HH_0$.
Indeed, $\HH_0$ is algebraically integrable because $\HH$ is,
and $\HH_0\subset\G$, so $\HH_0\subset\G^a$.
Conversely, $q^{-1}(\G^a)$ is an algebraically integrable
subfoliation of $\F$, and therefore
\[
q^{-1}(\G^a)\subset\F^a=\HH=q^{-1}\HH_0.
\]
It follows that $\G^a=\HH_0$.

The foliation $\HH_0$ is not uniruled.
This is immediate if $\HH_0=0$.
Otherwise, $\HH_0$ is regular, hence canonical, and
$K_{\HH_0}\equiv0$, so the assertion follows from
\cite[Corollary~3.8]{LPT18}.
Consequently, $\G$ is not uniruled either:
a covering family of rational curves tangent to $\G$ would be
tangent to its algebraic part $\G^a=\HH_0$, contradicting the
preceding conclusion.
Since $K_{\G}\equiv0$, another application of
\cite[Corollary~3.8]{LPT18} shows that $\G$ is canonical.

Lemma~\ref{lem:dominant} shows $p$-closedness for almost all primes
of $\G$.
Theorem~\ref{thm:numtriv} therefore makes $\G$, and hence $\F$,
algebraically integrable.

We have thus proved that $\F=\HH$.
In particular, $\F$ has a compact leaf.
Ou's compact-leaf corollary, included in Theorem~\ref{thm:ou},
proves regularity and an involutive splitting
\[
T_X=\F\oplus\E.
\]

\end{proof}
\begin{corollary}\label{cor:small-sing}
Let $X$ be a smooth complex projective variety of dimension at least two,
and let $\F\subset T_X$ be a corank one foliation.
Assume that $-K_{\F}$ is nef and that $\F$ is $p$-closed for almost all
primes. Suppose that one of the following conditions holds:
\begin{enumerate}
\item $\F$ has only simple singularities;
\item $\codim_X\Sing(\F)\geq3$.
\end{enumerate}
Then $\F$ is algebraically integrable and regular.
Moreover, there is a regular rank-one foliation $\E$ such that
\[
T_X=\F\oplus\E.
\]
\end{corollary}
\begin{proof}
We show that $\F$ admits a local holomorphic first integral at every point.
At a regular point this follows from the holomorphic Frobenius theorem.

In case (i),
since $\F$ is $p$-closed for almost all primes, it is in particular analytically
$p$-closed at every point in the sense of \cite[Definition~3.17]{Liu26}.
Thus \cite[Corollary~1.7]{Liu26} implies the existence of a holomorphic first integral
at each singular point.

In case (ii), choose a local holomorphic one-form $\omega$ defining
$\F$, with no divisorial zeroes.
Its zero locus has codimension at least three.
By \cite[Th\'eor\`eme~(0.1)]{Mal76}, one can write
$\omega=u\,dh$, where $u$ is a holomorphic unit and $h$ is a
nonconstant holomorphic function.
Hence $h$ is a local holomorphic first integral.

In either case, Theorem \ref{main:nef} (ii) applies and proves all the stated conclusions.
\end{proof}
\section{Foliated pairs on threefolds}\label{sec:foliated-pairs}
In this section, we give the proof of Theorem \ref{main:log-nef}. We first prove two lemmas comparing the canonical divisor of a foliation and its algebraic part. 

\begin{lemma}\label{lem:pair-algebraic-comparison}
Let $X$ be a normal $\QQ$-factorial projective threefold, and let $\F$
be a corank one foliation. Suppose that its algebraic part $\HH=\F^a$
has rank one and is induced by an almost holomorphic map. Let $B$ be an
effective $\QQ$-divisor with coefficients in $[0,1]$. If
$-(K_{\F}+B)$ is nef, then $K_{\HH}-K_{\F}$ is pseudo-effective.
\end{lemma}
\begin{proof}
Let $f:X\dashrightarrow T$ be the almost holomorphic quotient inducing
$\HH$, with connected general fibres, and let $\G_T$ be the quotient
foliation, so that $\F=f^{-1}\G_T$. Since $\dim T=2$, a general fibre
$C$ avoids $\operatorname{Sing}(X)$ and $\operatorname{Sing}(\F)$. So by adjunction, we have
\begin{equation}\label{eq:pair-leaf-restriction}
K_{\HH}|_C\sim K_C,
\qquad K_{\F}|_C\sim K_C.
\end{equation}

The general log leaf of $\HH$ is therefore $(C,0)$.

The horizontal components of $B$ are generically finite over $T$.
After shrinking the base further, their restrictions are smooth,
pairwise disjoint and \'etale over the base; the vertical components do
not meet the general fibre. Consequently, $B|_C$ has distinct support
points with coefficients at most one.

Put $P=-(K_{\F}+B)$, and choose an ample Cartier divisor $A$ on $X$.
For every positive rational number $\varepsilon$, the divisor
$P+\varepsilon A$ is ample. Choose a sufficiently divisible
integer $m>0$ such that $m(P+\varepsilon A)$ is very ample.
For a general member $H\in|m(P+\varepsilon A)|$, put
$D_\varepsilon=\frac1mH$.
On a general fibre $C$, the divisor $H|_C$ is reduced and its
support is disjoint from that of $B|_C$.
Thus all coefficients of $B|_C+D_\varepsilon|_C$ are at most one,
so $(C,B|_C+D_\varepsilon|_C)$ is log canonical. By \eqref{eq:pair-leaf-restriction},
\[
K_C+B|_C+D_\varepsilon|_C\equiv\varepsilon A|_C,
\]
which has positive degree. Since $X$ is $\QQ$-factorial, $K_{\HH}$
and $B+D_\varepsilon$ are $\QQ$-Cartier. We may therefore apply
\cite[Proposition~4.1]{Dru17} to $\HH$ and $B+D_\varepsilon$.
It follows that
\[
K_{\HH}+B+D_\varepsilon
\equiv K_{\HH}-K_{\F}+\varepsilon A
\]
is pseudo-effective. Letting $\varepsilon$ tend to zero proves the lemma.
\end{proof}

\begin{lemma}\label{lem:pair-nef-surface-quotient}
Let $X$ be a normal $\QQ$-factorial projective threefold, let $S$ be a
normal projective surface, and let $\varphi:X\dashrightarrow S$ be a
dominant rational map with connected general fibres. Let $\G$ be a
rank-one foliation on $S$ with canonical singularities and nef canonical
divisor. Set $\F=\varphi^{-1}\G$, and let $\HH$ be the foliation induced
by $\varphi$. If $K_{\HH}-K_{\F}$ is pseudo-effective, then
\(
K_{\HH}\equiv K_{\F}\)
and \(K_{\G}\equiv0.
\)
\end{lemma}
\begin{proof}
Resolve the rational map and apply flattening. After a further resolution
of the base and normalization of the main component, we obtain a
commutative diagram
\[
\begin{tikzcd}[column sep=large,row sep=large]
Z \arrow[r,"\nu"] \arrow[d,"q"'] &
X \arrow[d,dashed,"\varphi"]\\
S' \arrow[r,"\sigma"'] & S,
\end{tikzcd}
\]
where $Z$ is normal, $S'$ is smooth, $\nu$ and $\sigma$ are projective
birational morphisms, and $q$ is equidimensional with connected fibres.
Write
\[
\F_Z=\nu^{-1}\F=q^{-1}\G',\qquad
\HH_Z=\nu^{-1}\HH,\qquad \G'=\sigma^{-1}\G.
\]
Canonicality of $\G$ implies
\[
K_{\G'}=\sigma^*K_{\G}+E_S,
\qquad E_S\geq0.
\]
By \cite[Lemma~5.3]{CS25},
\[
K_{\F_Z}-K_{\HH_Z}
=q^*K_{\G'}+R(q)_{\mathrm{n\text{-}inv}},
\]
where $R(q)_{\mathrm{n\text{-}inv}}$ is the part of the ramification
divisor of $q$ whose components are not $\F_Z$-invariant.
In particular, this divisor is effective. Pushing forward to $X$ implies
\begin{equation}\label{eq:pair-quotient-comparison}
K_{\F}-K_{\HH}\sim_{\QQ}N+R,
\end{equation}
where
\[
N=\nu_*q^*\sigma^*K_{\G},
\qquad
R=\nu_*\bigl(q^*E_S+R(q)_{\mathrm{n\text{-}inv}}\bigr)\geq0.
\]
The divisor $N$ is pseudo-effective, since it is the birational
push-forward of a nef divisor. Thus $K_{\F}-K_{\HH}$ and its negative
are both pseudo-effective. The pseudo-effective cone contains no line,
so $K_{\F}-K_{\HH}\equiv0$. Applying the same observation to
\eqref{eq:pair-quotient-comparison} yields $N\equiv0$ and $R\equiv0$;
as $R$ is effective, $R=0$.

It remains to deduce $K_{\G}\equiv0$ from $N\equiv0$.
Take a smooth projective resolution $\tau:W\to Z$, and put
\[
\rho=\nu\circ\tau,
\qquad g=\sigma\circ q\circ\tau,
\qquad M=g^*K_{\G}.
\]
Then $M$ is nef and $\rho_*M=N$. Since $X$ is $\QQ$-factorial,
the divisor
\[
D=M-\rho^*N
\]
is $\QQ$-Cartier, $\rho$-exceptional and $\rho$-nef. The negativity
lemma \cite[Lemma~3.39]{KM98} shows $D\leq0$. Since $N\equiv0$,
we also have $M\equiv D$. For an ample divisor $A_W$ on $W$,
\[
0\leq M\cdot A_W^2=D\cdot A_W^2\leq0.
\]
A nonzero negative effective divisor has strictly negative intersection
with $A_W^2$, so $D=0$ and $g^*K_{\G}=M\equiv0$. Hence $K_{\G}\equiv0$.
\end{proof}

\begin{proof}[Proof of Theorem~\ref{main:log-nef}]
By \cite[Theorem~11.1]{CS21}, there is an F-terminalisation
$\mu:Y\to X$ such that $Y$ is $\QQ$-factorial and klt,
$\F_Y=\mu^{-1}\F$ has canonical singularities, and
\begin{equation}\label{eq:pair-terminalisation}
K_{\F_Y}+\mu_*^{-1}\Delta+E
=\mu^*(K_{\F}+\Delta),\qquad E\geq0,
\end{equation}
where $E$ is $\mu$-exceptional. Put $B=\mu_*^{-1}\Delta+E$.
The equality \eqref{eq:pair-terminalisation} shows that $(\F_Y,B)$
is log canonical, so the coefficients of $B$ belong to $[0,1]$.
Moreover,
\[
P:=-(K_{\F_Y}+B)=\mu^*\bigl(-(K_{\F}+\Delta)\bigr)
\]
is nef. By Lemma~\ref{lem:open}, $\F_Y$ is $p$-closed for almost
all primes. It is enough to prove that $\F_Y$ is algebraically integrable.

Suppose otherwise, and set $\HH=\F_Y^a$.
If $\HH=0$, apply \cite[Corollary~4.13]{CP19} on a resolution of $Y$
and push forward to see that $K_{\F_Y}$ is pseudo-effective.
Its negative is also pseudo-effective, since $-K_{\F_Y}=P+B$.
Hence $K_{\F_Y}\equiv0$, and Theorem~\ref{thm:numtriv} implies $\F_Y$
algebraically integrable, a contradiction. Since $\F_Y$ has rank two
and is not algebraically integrable, we obtain $\rank\HH=1$.

By \cite[Theorem~5.4]{CS25}, $\HH$ is induced by an almost holomorphic
map $Y\dashrightarrow S_0$ to a normal projective surface.
Lemma~\ref{lem:pair-algebraic-comparison} therefore implies $K_\HH-K_{\F_Y}$ is pseudo-effective.

Let $\G_0$ be the purely transcendental quotient foliation on $S_0$. By Seidenberg's theorem \cite{Sei68}, we may choose a smooth
projective model of $S_0$ on which the transformed foliation
has canonical singularities.
The transformed foliation remains purely transcendental, and hence
has pseudo-effective canonical divisor by
\cite[Corollary~4.13]{CP19}.
We then run the foliated surface MMP
\cite[Sections~III.1--III.2]{McQ08}.
Since the canonical divisor of the foliation is pseudo-effective,
the MMP terminates with a projective model $(S,\G)$ such that
$K_{\G}$ is nef.
Moreover, $S$ is klt and $\G$ has canonical singularities.
Let $\varphi:Y\dashrightarrow S$ be the induced rational map.
Then $\F_Y=\varphi^{-1}\G$, and $\HH$ is the foliation induced by
$\varphi$. Lemma~\ref{lem:pair-nef-surface-quotient} and
pseudo-effectivity of $K_\HH-K_{\F_Y}$ imply that $K_{\G}\equiv0$.

By Lemma~\ref{lem:dominant}, $\G$ is $p$-closed for almost all primes.
Since $S$ is klt and $\G$ is canonical, $\G$ is algebraically integrable by Theorem~\ref{thm:numtriv}, which contradicts pure
transcendence. Thus $\F_Y$, and hence $\F$, is algebraically integrable.
\end{proof}

\begin{corollary}\label{cor:lc-nef-threefold}
Let $X$ be a normal complex projective threefold, and let $\F$ be a
corank one foliation with log canonical singularities. If $-K_{\F}$
is nef and $\F$ is $p$-closed for almost all primes, then $\F$ is
algebraically integrable. In particular, this applies to log canonical
$\QQ$-Fano foliations on threefolds.
\end{corollary}
\begin{proof}
Apply Theorem~\ref{main:log-nef} with $\Delta=0$.
\end{proof}

\section{Applications}

We now give some applications of the preceding results.
We first prove Conjecture~\ref{conj:esbt} for certain classes of
Fano foliations, starting with foliations of large index in
Corollary~\ref{cor:large-index} and then turning to Fano foliations
on fourfolds with Picard number one in Corollary~\ref{cor:fourfold}.
In preparation for the latter, we prove the conjecture for
log canonical degree-two foliations on $\PP^n$ and for log canonical
Mukai foliations on varieties of Picard number one; see
Corollaries~\ref{cor:degree-two} and~\ref{cor:mukai}, respectively.
Finally, we prove the adjoint Fano criterion stated in
Corollary~\ref{cor:intro-adjoint}.
\subsection{Fano foliations of large index}
Here the ambient variety may be singular.
A $\QQ$-Fano foliation means that $K_{\F}$ is $\QQ$-Cartier and $-K_{\F}$
is ample. Its index is
\begin{equation}\label{eq:index}
\iota(\F)=\max\{a\in\QQ_{>0}:-K_{\F}\smq aH
\text{ for an ample Cartier divisor }H\}.
\end{equation}

\begin{theorem}[\cite{Liu23}, Theorem 1.3]\label{thm:liu-index}
Let $X$ be a $\QQ$-factorial normal projective variety and
$\F\subsetneq T_X$ a $\QQ$-Fano foliation. Then $r^a(\F)\geq\iota(\F)$.
If equality holds, $X$ is a normal generalized cone over a polarized
$\QQ$-factorial projective variety $(T,L)$, with vertex $\PP^{r^a(\F)-1}$,
and $\F$ is the pull-back of a foliation $\G$ on $T$ satisfying $K_{\G}\smq0$.
\end{theorem}

In the equality case of Theorem~\ref{thm:liu-index}, put
$k=r^a(\F)$. The generalized cone structure defines a birational
morphism
\begin{equation}\label{eq:cone}
\begin{tikzcd}[column sep=large,row sep=large]
\PP_T(L\oplus\OO_T^{\oplus k})\arrow[r,"\mu"]\arrow[d,"q"']&
X\arrow[dl,dashed,"h"]\\
T&
\end{tikzcd}
\end{equation}
where $q$ is the natural projection and $h=q\circ\mu^{-1}$.
The morphism $\mu$ contracts
$E=\PP_T(\OO_T^{\oplus k})\simeq T\times\PP^{k-1}$
onto the vertex $\PP^{k-1}$ by the second projection,
and is an isomorphism outside $E$.
Thus the complement of the vertex in $X$ is an affine-space
bundle of relative dimension $k$ over $T$;
see \cite[Example~4.4]{Liu23}.
\begin{lemma}\label{lem:cone-klt}
If a normal generalized cone $X$ is klt, then its base $T$ is klt.
\end{lemma}
\begin{proof}
By the description above, the open complement of the vertex in $X$ is
an affine-space bundle over $T$. Zariski locally on $T$, it is isomorphic to
$U\times\AA^k$. The klt property is invariant under smooth surjective local
products. Concretely, pulling a resolution of $U$ back to $U\times\AA^k$
preserves every discrepancy and introduces no change in its coefficient.
Since $X$ is klt, all these discrepancies are greater than $-1$, and the
same is true on $U$. Thus $T$ is klt.
\end{proof}
\begin{corollary}\label{cor:large-index}
Let $X$ be a $\QQ$-factorial klt projective variety of dimension $n\geq3$,
and let $\F$ be a corank one $\QQ$-Fano foliation.
If $\F$ is $p$-closed for almost all primes and $\iota(\F)\geq n-2$,
then $\F$ is algebraically integrable.
\end{corollary}
\begin{proof}
Assume that $\F$ is not algebraically integrable.
Since its rank is $n-1$, one has $r^a(\F)\leq n-2$.
The index bound and Theorem~\ref{thm:liu-index} therefore imply
\begin{equation}\label{eq:index-chain}
n-2\geq r^a(\F)\geq\iota(\F)\geq n-2.
\end{equation}
All terms are equal. We are in the generalized-cone case of
Theorem~\ref{thm:liu-index}, and its base $T$ has dimension two.
Let $\G$ be the resulting foliation on $T$. It has rank one and $K_{\G}\smq0$.

The foliation induced by $h$ is an algebraically integrable
subfoliation of $\F$ of rank $n-2=r^a(\F)$.
It therefore coincides with $\F^a$, so the quotient foliation
$\G$ is purely transcendental.
Lemma~\ref{lem:cone-klt} shows that $T$ is klt.
By Lemma~\ref{lem:dominant}, $\G$ is $p$-closed for almost all primes, which is a contradiction with Corollary~\ref{cor:pure-numtriv}. So $\F$ is algebraically integrable.
\end{proof}
\begin{corollary}
Let $X$ be a $\QQ$-factorial klt projective threefold and let $\F$ be a corank
one $\QQ$-Fano foliation. If $\F$ is $p$-closed for almost all primes and
$\iota(\F)\geq1$, then it is algebraically integrable.
In particular this holds whenever $K_{\F}$ is Cartier.
\end{corollary}
\begin{proof}
Apply Corollary~\ref{cor:large-index} with $n=3$.
If $K_{\F}$ is Cartier, $H=-K_{\F}$ is an ample Cartier divisor
and the value $a=1$ is allowed in \eqref{eq:index}.
\end{proof}

\subsection{Degree-two foliations on projective space}
We follow the same definition as \cite[Section~9.1]{MP24}. For a codimension one foliation $\F$ of degree $d$ on $\PP^n$, the normal
bundle and canonical class are
\begin{equation}\label{eq:degree}
N_{\F}\simeq\OO_{\PP^n}(d+2),\qquad
K_{\F}\simeq\OO_{\PP^n}(d-n+1).
\end{equation}
Let $\pi:\AA^{n+1}\setminus\{0\}\to\PP^n$ be the natural projection.
By the Euler sequence, the defining section
$\omega_{\F}\in H^0(\PP^n,\Omega_{\PP^n}^1(d+2))$
is represented by a homogeneous one-form
\[
\omega=\sum_{i=0}^{n}A_i\,dx_i
\]
defining $\pi^{-1}\F$, where the $A_i$ are homogeneous polynomials
of degree $d+1$ satisfying
\[
\sum_{i=0}^{n}x_iA_i=0.
\]
The coefficients $A_i$ have no common nonconstant factor because
$\omega_{\F}$ has no divisorial zeroes.
See \cite[Sections~9.1--9.2]{MP24}.
\begin{lemma}\label{lem:projection}
Let $n>m\geq2$, let $\pi:\PP^n\dashrightarrow\PP^m$ be a linear projection,
and let $\G$ be a corank one foliation of degree $d$ on $\PP^m$.
Set $\F=\pi^{-1}\G$. If $\sigma:W\to\PP^n$ is the blow-up of the centre
of $\pi$, with exceptional divisor $E$, then $E$ is non-invariant for
$\sigma^{-1}\F$ and
\begin{equation}\label{eq:projection-discrepancy}
a(E,\F)=m-d-2.
\end{equation}
Consequently, log canonicity of $\F$ implies $d\leq m-1$.
\end{lemma}
\begin{proof}
Choose coordinates so that the projection is given by
$[x_0:\cdots:x_n]\mapsto[x_0:\cdots:x_m]$ and its centre is
$\Lambda=(x_0=\cdots=x_m=0)$.
Write a saturated defining form of $\G$ as
\[
\omega=\sum_{i=0}^m A_i(x_0,\ldots,x_m)\,dx_i,
\quad \deg A_i=d+1,\quad \sum_{i=0}^m x_iA_i=0.
\]
The same form defines $\F$ on $\PP^n$; extending the coefficient ring by
the remaining variables introduces no common divisorial factor.

Work near a general point of $\Lambda$ and in the blow-up chart
$x_0=u$, $x_i=uv_i$ for $1\leq i\leq m$.
The remaining coordinates are coordinates along the centre. Then
\[
\sigma^*\omega
=u^{d+1}\left(A_0(1,v)+\sum_{i=1}^m v_iA_i(1,v)\right)du
+u^{d+2}\sum_{i=1}^m A_i(1,v)\,dv_i.
\]
The Euler relation makes the first term vanish. The second term has exact
order $d+2$ along $u=0$: its residual form is not identically zero at the
generic point of $E$. Saturating the pull-back therefore removes the
factor $u^{d+2}$. So
\[
N_{\sigma^{-1}\F}\simeq\sigma^*N_{\F}\otimes\OO_W(-(d+2)E).
\]
The centre has codimension $m+1$, so $K_W=\sigma^*K_{\PP^n}+mE$.
So $a(E,\F)=m-d-2$ which proves \eqref{eq:projection-discrepancy}.

The saturated form $\sum A_i(1,v)\,dv_i$ annihilates $\partial_u$,
which is transverse to $E$. Thus $E$ is non-invariant.
Log canonicity requires $m-d-2\geq-1$, proving the last assertion.
\end{proof}
\begin{corollary}\label{cor:projection-lc}
A linear pull-back to $\PP^n$, $n\geq3$, of a degree-two foliation on
$\PP^2$ is not log canonical.
\end{corollary}
\begin{proof}
In Lemma~\ref{lem:projection}, take $m=d=2$. The non-invariant divisor $E$
has discrepancy $-2$.
\end{proof}

We will use the following consequence of the complete classification of Cerveau--Lins Neto.
\begin{theorem}[\cite{CLN96}, Theorem A]\label{thm:degree-two-class}
A codimension one foliation of degree two on $\PP^n$, $n\geq3$, is either
defined by a closed rational one-form or is the linear pull-back of a
degree-two foliation on $\PP^2$.
\end{theorem}

\begin{corollary}\label{cor:degree-two}
Let $\F$ be a corank one foliation of degree two on $\PP^n$, where $n\geq3$.
If $\F$ has log canonical singularities and is $p$-closed for almost all primes,
then $\F$ is algebraically integrable.
\end{corollary}
\begin{proof}
By Corollary~\ref{cor:projection-lc} and Theorem~\ref{thm:degree-two-class}, $\F$ is defined by a closed rational one-form on $\PP^n$. So algebraic integrability follows from Theorem~\ref{thm:closed}.
\end{proof}

\subsection{Mukai foliations on varieties of Picard number one}
A Fano foliation of rank $r\geq3$ is called a Mukai foliation if its index
is $r-2$. Thus a corank one Mukai foliation on an $n$-fold has index $n-3$,
with $n\geq4$. For a corank one foliation on $\PP^n$, $n\geq4$,
the Mukai condition is equivalent to having degree two
by \eqref{eq:degree}.

We will use the following classification by Araujo and Druel.
\begin{theorem}[\cite{AD17}, Theorem 5]\label{thm:mukai-class}
Let $X$ be a smooth projective $n$-fold with $n\geq4$, $\rho(X)=1$,
and $X\not\simeq\PP^n$. If $X$ carries a corank one Mukai foliation $\F$,
then one of the following holds:
\begin{enumerate}
\item $X\simeq Q^n$, and $\F$ is induced by a pencil of ambient quadrics
containing a double hyperplane;
\item $X\simeq Q^n$, and $\F$ is the pull-back, under the restriction of a
linear projection $\PP^{n+1}\dashrightarrow\PP^2$, of a rank-one foliation
induced by a global vector field on $\PP^2$;
\item $X$ is Fano of index $n-1$, and $\F$ is induced by a pencil in the ample
generator of $\Pic(X)$.
\end{enumerate}
\end{theorem}
The first and third alternatives are already algebraically integrable.
We only need to consider the second case.
\begin{lemma}\label{lem:vectorfield}
Let $\G$ be a rank-one foliation on $\PP^2$ induced by a nonzero global
vector field. If $\G$ is $p$-closed for almost all primes,
then it is algebraically integrable.
\end{lemma}
\begin{proof}
Let $B\geq0$ be the divisorial part of the zero locus of the vector field.
The saturated subsheaf of $T_{\PP^2}$ generated by the vector field
is $\G\simeq\OO_{\PP^2}(B)$, so \eqref{eq:degree} implies
\[
1-d=\deg B\geq0,
\]
where $d$ is the degree of $\G$.
Thus $d\leq1$. If $d=0$, then $\G$ is induced by a pencil of lines
through a point and is therefore algebraically integrable.
If $d=1$, its $p$-closedness for almost all primes implies algebraic
integrability by \cite[Lemma~11.5]{MP24}.
\end{proof}
\begin{corollary}\label{cor:mukai}
Let $X$ be a smooth projective variety of dimension $n\geq4$ with $\rho(X)=1$,
and let $\F$ be a corank one Mukai foliation.
Assume that $\F$ has log canonical singularities and is $p$-closed for
almost all primes. Then $\F$ is algebraically integrable.
\end{corollary}
\begin{proof}
If $X\simeq\PP^n$, then $\F$ has degree two, so its algebraic
integrability follows from Corollary~\ref{cor:degree-two}.
We may therefore assume that $X\not\simeq\PP^n$.

By Theorem~\ref{thm:mukai-class}, either $\F$ is induced by a pencil,
in which case it is algebraically integrable, or
$\F=\pi^{-1}\G$, where $\pi:X\dashrightarrow\PP^2$ is a dominant
rational map and $\G$ is induced by a nonzero global vector field
on $\PP^2$.
In the latter case, Lemma~\ref{lem:dominant} implies that $\G$
is $p$-closed for almost all primes.
Lemma~\ref{lem:vectorfield} then shows that $\G$ is algebraically
integrable, and hence so is $\F$.
\end{proof}
\begin{remark}
In Corollary~\ref{cor:mukai}, log canonicity is used only when $X\simeq\PP^n$.
For $X\not\simeq\PP^n$, the same proof shows algebraic integrability without
a singularity assumption on $\F$.
\end{remark}
\begin{proof}[Proof of Corollary~\ref{cor:intro-degree-mukai}]
Parts (i) and (ii) are Corollaries~\ref{cor:degree-two}
and~\ref{cor:mukai}, respectively.
\end{proof}
In dimension four, every positive integral Fano index is either at least
two or is the Mukai index. We obtain the following consequence.
\begin{corollary}\label{cor:fourfold}
Let $X$ be a smooth projective fourfold with $\rho(X)=1$, and let $\F$ be a
corank one Fano foliation with log canonical singularities.
If $\F$ is $p$-closed for almost all primes, then it is algebraically integrable.
\end{corollary}
\begin{proof}
Since $X$ is smooth, $-K_{\F}$ is Cartier, and hence its Fano index
$a=\iota(\F)$ is a positive integer.
If $a\geq2=\dim X-2$, then Corollary~\ref{cor:large-index}
implies that $\F$ is algebraically integrable.

It remains to consider $a=1$.
Since $\F$ has rank three, $\F$ is a Mukai foliation.
Its algebraic integrability therefore follows from
Corollary~\ref{cor:mukai}.
\end{proof}
\begin{proof}[Proof of Corollary~\ref{cor:intro-fano}]
Part (i) is Corollary~\ref{cor:large-index}.
For part (ii), the threefold case follows from
Theorem~\ref{main:log-nef} with $\Delta=0$,
since $-K_{\F}$ is ample and hence nef.
The fourfold case is Corollary~\ref{cor:fourfold}.
\end{proof}

\subsection{Adjoint Fano foliated structures}
For a foliation $\F$ on a smooth projective variety $X$ and
$t\in[0,1]\cap\QQ$, set
\[
K_t=tK_{\F}+(1-t)K_X.
\]
We consider the Fano condition that $-K_t$ is ample for the
adjoint foliated structure $(X,\F,t)$; see
\cite[Section~1 and Definition~3.6]{CHL25}.
We first show that, for $t$ sufficiently close to one, this condition
implies that $-K_{\F}$ is nef.
We then apply Theorem~B to prove the adjoint Fano criterion
stated in Corollary~\ref{cor:intro-adjoint}.
\begin{lemma}\label{lem:convex}
Let $\mathcal C$ be a closed convex cone in a finite-dimensional real vector
space, and let $h$ lie in the interior of its dual cone, with
$h(\gamma)>0$ for $\gamma\in\mathcal C\setminus\{0\}$.
Then $\Sigma=\{\gamma\in\mathcal C:h(\gamma)=1\}$ is compact.
Each extreme point $\beta$ of $\Sigma$ spans an extremal ray of $\mathcal C$.
\end{lemma}
\begin{proof}
If $\mathcal C=\{0\}$, both assertions are immediate.
Otherwise, fix a norm on the ambient finite-dimensional vector space.
Since $\mathcal C$ is closed, the set
\[
S=\{\gamma\in\mathcal C:\|\gamma\|=1\}
\]
is compact. The linear functional $h$ is strictly positive on $S$,
so
\[
\varepsilon:=\min_{\gamma\in S}h(\gamma)>0.
\]
For every $\gamma\in\Sigma$, we have $h(\gamma)=1$ and
$\gamma/\|\gamma\|\in S$. Hence
\[
1=h(\gamma)
=\|\gamma\|\,h\!\left(\frac{\gamma}{\|\gamma\|}\right)
\geq\varepsilon\|\gamma\|.
\]
Thus $\Sigma$ is bounded. It is also closed, being the intersection
of $\mathcal C$ with the hyperplane $\{h=1\}$, and is therefore compact.

Now let $\beta$ be an extreme point of $\Sigma$.
To show that $\RR_{\geq0}\beta$ is an extremal ray of $\mathcal C$,
it suffices, after rescaling, to consider a decomposition
$\beta=u+v$ with $u,v\in\mathcal C$ and prove that both vectors
belong to $\RR_{\geq0}\beta$.
This is immediate if either vector is zero.
Otherwise, set $a=h(u)>0$ and $b=h(v)>0$.
Since $h(\beta)=1$, we have $a+b=1$, and
\[
\beta=a\left(\frac{u}{a}\right)
      +b\left(\frac{v}{b}\right),
\qquad
\frac{u}{a},\frac{v}{b}\in\Sigma.
\]
This expresses $\beta$ as a convex combination of two points of
$\Sigma$ with positive coefficients.
Since $\beta$ is an extreme point, both points must equal $\beta$.
Therefore $u=a\beta$ and $v=b\beta$, as required.
\end{proof}
\begin{proposition}\label{prop:threshold}
Let $X$ be a normal projective klt variety of dimension $n$, let $M$ be a
$\QQ$-Cartier divisor, and let $r>0$ be an integer such that $rM$ is Cartier.
Suppose that $-(tM+(1-t)K_X)$ is ample for some rational $t$ satisfying
\begin{equation}\label{eq:klt-threshold}
\frac{2nr}{2nr+1}\leq t<1.
\end{equation}
Then $-M$ is nef. If $X$ is smooth, the same conclusion holds under the
weaker bound
\begin{equation}\label{eq:smooth-threshold}
\frac{(n+1)r}{(n+1)r+1}\leq t<1.
\end{equation}
\end{proposition}
\begin{proof}
Set $H=-(tM+(1-t)K_X)$, which is ample by assumption.
Suppose for contradiction that $-M$ is not nef.
Then $M$ has positive intersection with some curve, so its maximum
on the compact slice
\[
\Sigma=\{\gamma\in\NE(X):H\cdot\gamma=1\}
\]
is positive. Denote this maximum by $\lambda>0$.
The subset
\[
\Sigma_{\max}
=\{\gamma\in\Sigma:M\cdot\gamma=\lambda\}
\]
is a nonempty compact face of $\Sigma$.
Choose an extreme point $\beta$ of $\Sigma_{\max}$.
Since $\Sigma_{\max}$ is a face, $\beta$ is also an extreme point
of $\Sigma$. Lemma~\ref{lem:convex} therefore shows that
\[
R=\RR_{\geq0}\beta
\]
is an extremal ray of $\NE(X)$, and $M\cdot\beta=\lambda>0$.

Since $H\cdot\beta=1$, we have
\[
(1-t)K_X\cdot\beta=-1-t\lambda<0.
\]
Thus $R$ is $K_X$-negative.
The cone theorem implies the existence of a rational curve $C$ spanning $R$ such that
\[
0<-K_X\cdot C\leq b.
\]
Here we take $b=2n$ when $X$ is klt, by
\cite[Theorem~3.7]{KM98} and \cite[Theorem~1]{Kaw91},
and $b=n+1$ when $X$ is smooth, by
\cite[Theorem~1.1]{Fuj25}.

Because $rM$ is Cartier and $M\cdot C>0$, the intersection
$rM\cdot C$ is a positive integer. Hence $M\cdot C\geq1/r$.
In either case, the assumed bound on $t$ is
$t\geq br/(br+1)$, so
\[
\begin{aligned}
-H\cdot C
&=tM\cdot C+(1-t)K_X\cdot C\\
&\geq \frac{t}{r}-b(1-t)
\geq0.
\end{aligned}
\]
This contradicts $H\cdot C>0$, since $H$ is ample.
Therefore $-M$ is nef.
\end{proof}

\begin{lemma}\label{lem:limit}
Let $X$ be a projective variety, and let $M$ and $K_X$ be
$\QQ$-Cartier divisors.
The following conditions are equivalent:
\begin{enumerate}
\item $-(sM+(1-s)K_X)$ is ample for every rational
$s\in(0,1)$ sufficiently close to one;
\item $-M$ is nef and $-(tM+(1-t)K_X)$ is ample for
at least one rational $t\in(0,1)$.
\end{enumerate}
\end{lemma}
\begin{proof}
Assume (i). As $s\to1$ from below, the ample classes
$-(sM+(1-s)K_X)$ converge to $-M$, so $-M$ is nef.

Conversely, assume (ii) and fix such a rational number $t$.
For every rational $s$ with $t\leq s<1$, we have
\begin{equation}\label{eq:convex-combination}
-\bigl(sM+(1-s)K_X\bigr)
=\frac{1-s}{1-t}\bigl(-tM-(1-t)K_X\bigr)
+\frac{s-t}{1-t}(-M).
\end{equation}
The first summand is ample, since its coefficient is positive,
and the second is nef.
Their sum is therefore ample, proving (i).
\end{proof}

\begin{proof}[Proof of Corollary \ref{cor:intro-adjoint}]
Since $X$ is smooth, $K_{\F}$ is Cartier.
The smooth case of Proposition~\ref{prop:threshold}, applied with
$M=K_{\F}$ and $r=1$, shows that $-K_{\F}$ is nef.
The conclusion follows from Theorem \ref{main:nef}.
\end{proof}
\begin{remark}
By Proposition~\ref{prop:threshold} and Lemma~\ref{lem:limit},
the ampleness hypothesis in Corollary~\ref{cor:intro-adjoint} is equivalent
to requiring $-K_s$ to be ample for every rational $s\in(0,1)$
sufficiently close to one.    
\end{remark}

\bibliographystyle{amsalpha}
\bibliography{references}

\end{document}